\documentstyle[12pt]{article} 

\def\ds{\displaystyle}
\def\be{\begin{equation}} 
\def\ee{\end{equation}} 
\newtheorem{th}{Theorem}[section]  
\newtheorem{lm}[th]{Lemma}
\newtheorem{cor}[th]{Corollary}
\newtheorem{pr}[th]{Proposition}

\def\RR{{\hbox{I\kern-.2em\hbox{R}}}}
\def\PP{{\hbox{I\kern-.2em\hbox{P}}}}
\def\EE{{\hbox{I\kern-.2em\hbox{E}}}}
\def\ZZ{{\hbox{Z\kern-.4em\hbox{Z}}}}
\def\one{{\bf 1}}

\def\R{{\bf R}}

\def\C{{\bf C}}

\def\wt{\widetilde}

\def\qed{\hfill\hbox{\rule{6pt}{6pt}}}
\def\lg{\langle}
\def\rg{\rangle}

\def\cE{{\cal E}}
\def\cF{{\cal F}}
\def\cM{{\cal M}}
\def\cS{{\cal S}}

\def\one{{\bf 1}}

\begin{document} 
\setlength{\topmargin}{-4mm}
%
%
\begin{center}
{\LARGE \bf 
The sector constants of continuous state \\  
\vspace{3mm} 
branching processes with immigration }
\footnote{
Supported in part by 
JSPS, Grant-in-Aid for Scientific Research (A) No. 21244009. 
}
\end{center} 

%
%


\begin{center}
Kenji Handa 
\end{center} 

\begin{center}
{ 
Department of Mathematics \\
Saga University \\ 
Saga 840-8502 \\ 
Japan
} \\  
e-mail: 
handa@ms.saga-u.ac.jp \\ 
FAX: +81-952-28-8501
\end{center} 
\begin{center}
{\em Dedicated to Professor Tadahisa Funaki 
on the occasion of his 60th birthday}
\end{center}

%
%

\begin{center} 
\begin{minipage}[t]{12cm} 
\small 
Continuous state branching processes with immigration 
are studied. We are particularly concerned with 
the associated (non-symmetric) Dirichlet form. 
After observing that gamma distributions are 
only reversible distributions for this class of models, 
we prove that every generalized gamma convolution 
is a stationary distribution of the process
with suitably chosen branching mechanism 
and with continuous immigration. 
For such non-reversible processes, 
the strong sector condition is discussed 
in terms of a characteristic called the Thorin measure. 
In addition, some connections with notion from 
noncommutative probability theory will be pointed out 
through calculations involving the Stieltjes transform. 
\end{minipage}
\end{center}

%
%

\section{Introduction} 
\setcounter{equation}{0} 

Besides its significance in the physical context, 
the (time-)reversibility can be thought of as a mathematical condition 
which guarantees a certain kind of `solvability' of the equilibrium state 
and usually makes one possible to deduce explicit consequences. 
On the other hand, it is likely that 
the reversibility is a restrictive condition, and 
it fails for a number of stochastic models 
with stationary distributions of interest. 
In this paper, our attempt will be made 
in quantitative discussions on 
the degree of irreversibility of such systems. 
Let us illustrate roughly in a general setting 
the situation we will be concerned with.  
Suppose that we are given a Markov process 
with a stationary distribution $\nu$ 
and generator $L$, say.  
Then consider a bilinear form $\cE$ defined by 
\be 
\cE(f,g)= -\int Lf(x)g(x)\nu(dx),  
\qquad f,g\in D(L),  
                                      \label{1.1} 
\ee 
where $D(L)$ is the domain of $L$. 
The Dirichlet form is a suitable extension of $\cE$ and 
it is well-known that 
the symmetry of $\cE$ is interpreted as 
the reversibility of the Markov process. 
We say that $\cE$ satisfies the strong sector condition 
if there exists a finite constant $C$ such that 
\be  
\cE(f,g)\le C \cE(f,f)^{1/2}\cE(g,g)^{1/2}, 
\qquad f,g\in D(L).  
                                      \label{1.2} 
\ee 
Based on a weaker version called the weak sector condition, 
the theory of symmetric Dirichlet forms has been 
successfully extended to non-symmetric cases in \cite{MR92}. 
The strong sector condition is known also to play an essential role 
in the proof of the invariance principle for 
additive functionals of non-symmetric Markov processes. 
See \cite{OS95} and \cite{V95}. 
(See also \cite{KO03} for a generalization.) 
Intuitively, the validity of this condition tells us that 
the process is a small perturbation from a symmetric one. 
It seems typical that verification of (\ref{1.2}) 
depends heavily on the mathematical structure of the process. 
In our subsequent discussions it will be convenient to denote by 
${\rm Sect}(\cE)$ the infimum of $C$'s satisfying (\ref{1.2}) 
if any, and set ${\rm Sect}(\cE)=\infty$ otherwise. 
We call ${\rm Sect}(\cE)$ the sector constant of $\cE$. 
Clearly ${\rm Sect}(\cE)\ge 1$. 
If $\cE$ is symmetric, we have ${\rm Sect}(\cE)=1$. 
As expected, the converse holds true in general. 
(See Proposition 3.1 below for the proof.) 
Hence the difference ${\rm Sect}(\cE)-1$ can be thought of as 
a `degree' of asymmetry of $\cE$ and 
of irreversibility of the process. 
Our objective is to show 
such a property of ${\rm Sect}(\cE)$ 
in an explicit way for some specific class of models. 
We thus seek for an upper bound of ${\rm Sect}(\cE)-1$ 
which should be given in terms of 
certain characteristics of the models, 
and the bound is then required to 
vanish precisely in the reversible case. 

As a model which will be discussed in the present paper, 
we adopt the continuous state branching process with immigration 
(called also the CBI-process), 
which is a Markov process on $\R_+:=[0,\infty)$. 
Fundamental results of this process, 
including limit theorems from Galton-Watson processes with immigration 
and the complete determination of the generator, 
are obtained by Kawazu and Watanabe \cite{KW71}. 
Since then, this model has been studied extensively not only 
because of the rich and nice mathematical structure 
which has allowed us to obtain a number of concrete results of interest 
but also of its importance in various applications. 
The aforementioned authors showed its interesting applications 
in the context of stochastic analysis as well. 
In addition, since the CIR model \cite{CIR85}, 
a mathematical finance model for evolution of interest rate, 
is included as a special case (in fact, the diffusion case), 
the class of CBI-processes serves also as 
a useful generalization of the CIR model in such a context \cite{DFS}. 
We intend to reveal further aspects of the CBI-process 
regarding the non-reversible stationary distribution 
and the non-symmetric Dirichlet form. 

The time evolution of the CBI-process in general incorporates 
two kinds of dynamics; 
the one describing the branching of particles 
and the other being due to immigration. 
It is of essential importance to take 
into consideration the effect of immigration. 
One of consequences of introducing immigration is 
ergodicity of the process; 
it may exhibit the strong convergence to 
a unique stationary distribution, if any, 
as time goes to infinity. 
Actually, under suitable assumptions, 
the positivity of the spectral gap of $L$ follows 
from the result in \cite{Og70}. 
(See the discussion in the paragraph 
preceding to Lemma 2.1 below for the precise statement.) 
Furthermore, our process has so nice structure as 
to make it possible to get information 
of the stationary distribution through 
an explicit representation of the Laplace transform. 
This formula, a key tool throughout this paper, 
is due to Ogura \cite{Og70}, 
who carried out detailed calculations of 
the spectral representation for the CBI-process. However, 
it is typically difficult to deduce direct expressions 
(e.g. the density function) of the stationary distribution 
and so one needs to exploit other structures. 
(Among exceptions are gamma distributions, 
which are reversible distributions of the CIR models.) 
Another feature of the model which is crucial to us 
is the branching property, meaning that 
the law of the sum of two identical and independent 
processes starting from $x_1$ and $x_2$ respectively 
coincides with the law of a process starting from $x_1+x_2$. 
By virtue of this property, 
the law of the process at an arbitrarily fixed time 
is necessarily infinitely divisible and 
so is the stationary distribution. 
(See \cite{KRM} for recent studies of stationary 
distributions of the CBI-processes.) 
In \cite{Stannat05}, 
the Poincar\'e inequality 
for a class of infinitely divisible distributions 
on $\R_+$ (and more general spaces) was proved 
by reducing it to an analogous estimate 
for the associated L\'evy measure. 
It will turn out that a suitably modified 
argument works well for the sector constant estimate. 

Non-reversible stationary distributions 
we will focus attention on are generalized 
gamma convolutions \cite{Bondesson} (GGC's for short), 
namely weak limits of finite convolutions of gamma distributions. 
(See also Section 5, Chapter VI of \cite{SvH04} for general accounts 
and \cite{JRY08} for a recent survey and related topics.) 
We regard these distributions as `perturbations' 
from gamma distributions. 
A GGC without `translation term' is determined uniquely 
by the so-called Thorin measure, 
which appears in the logarithm of the Laplace transform 
and prescribes the weight of convolutions. 
For example, every gamma distribution has a degenerate Thorin measure. 
Therefore, the actual problems we are going to consider 
in the subsequent sections are outlined as follows.  \\ 
(I) Show that there does not exist 
a (nondegenerate) reversible distribution of 
the CBI-processes except gamma distributions. \\ 
(II) Given a GGC with Thorin measure $m$, 
choose a branching mechanism so that 
the CBI-process has the GGC 
as a unique stationary distribution. \\ 
(III) For the bilinear form $\cE$ associated with that process, 
give an upper bound $C=C(m)$ of ${\rm Sect}(\cE)$ 
such that $C(m)=1$ if and only if $m$ is degenerate. \\ 
We will see that the reversibility problem (I) reduces to solving 
certain functional equations involving `characteristics' 
of the mechanisms of branching and immigration. 
Our solution to (II) will turn out to rely on 
the theory of Bernstein functions \cite{SSV}. 
We also make use of Stieltjes transforms in order to get 
further information 
(e.g., the one needed to solve (III)) 
on the branching mechanism chosen. 
In this context some connections with 
notion from non-commutative probability theory 
(such as the so-called Boolean convolution 
and the free Poisson distribution) will be pointed out. 

The organization of this paper is as follows. 
In the next section a precise description of 
CBI-processes is given and then the problem (I) is solved. 
In Section 3, we present some basic results 
on the strong sector condition for a subclass of CBI-processes  
for which an integration by parts formula is available. 
In Section 4, both the problems (II) and (III) are solved by 
constructing the CBI-process associated with a GGC and 
then applying the results in Section 3. 
In Section 5, we give some examples to illustrate 
consequences of our results and discuss related topics. 

\section{The model and its stationary distribution}
\setcounter{equation}{0} 

Following \cite{KW71}, 
we begin with a precise description of our model, 
namely the CBI-process 
in terms of the generator. 
For the purpose of this paper, 
we shall restrict the discussion to 
a class of conservative CBI-processes. 
In view of Theorem 1.1$'$, Theorem 1.2 in \cite{KW71}, 
and results (especially Proposition 1.1) in \cite{Og70}, 
the assumptions made below are not optimal 
but useful in order that the results of this section 
are not more complicated than are necessary in the subsequent sections. 
(Recently, the detailed analysis of stationary distributions 
was done in \cite{KRM} for conservative CBI-processes.) 
The generator $L$ of our process takes the following form: 
\begin{eqnarray}
Lf(x)
&=& axf''(x)-bxf'(x) 
+x\int_0^{\infty}\left[f(x+y)-f(x)-yf'(x)\right]n_1(dy) 
\nonumber  \\ 
& & +\delta f'(x)+\int_0^{\infty}\left[f(x+y)-f(x)\right]n_2(dy), 
\qquad x\in\R_{+}, 
                                             \label{2.1} 
\end{eqnarray}
where $a\ge 0, b\ge 0, \delta\ge 0$, and 
measures $n_1$ and $n_2$ on $(0,\infty)$
are supposed to satisfy 
\be 
\int_0^{\infty}\min\{y^2,y\}n_1(dy)
+\int_{(0,1)}yn_2(dy)
+\int_{[1,\infty)}(1+\log y)n_2(dy)<\infty. 
                                             \label{2.2} 
\ee 
This process approximates 
(asymptotically critical) Galton-Watson branching processes 
with immigration in large population limit. 
In this context dynamical meaning of 
the constants and measures appearing in (\ref{2.1}) 
may be explained as follows. 
While $a$ is the asymptotic variance of the offspring distributions 
associated with the branching mechanisms, $b$ comes from 
the first order approximation to mean 1. 
$\delta$ is the rate of change in mean of immigrating population. 
$n_1$ and $n_2$ describe effects of 
big changes in population size which occur 
in `macroscopic time scale' and 
are caused by branch(-death) and immigration, respectively. 
To avoid the triviality in discussing the equilibrium of the model, 
we make the assumption implying that 
both branching and immigration mechanisms are actually present. 
To be precise, defining for $\lambda\ge 0$ 
\[ 
R(\lambda)= -a\lambda^2-b\lambda
-\int_0^{\infty}\left(e^{-\lambda y}-1+\lambda y\right)n_1(dy) 
\] 
and 
\[ 
F(\lambda)= \delta\lambda
 +\int_0^{\infty}(1-e^{-\lambda y})n_2(dy),  
\] 
we assume throughout that 
\be 
\mbox{neither}\  R\equiv 0 \ \mbox{nor}\  F\equiv 0. 
                                             \label{2.3} 
\ee 
The functions $R$ and $F$ are called 
the branching mechanism and the immigration mechanism, respectively, 
and their interplay will be crucial in the ergodic behavior of 
the CBI-process. 

As in the literature on CBI-processes, 
a large amount of calculations below will be based on the Laplace transforms, 
which can be expressed in terms of the associated 
$\Psi$-semigroup, a one-parameter family 
$\{\psi(t,\cdot)\}_{t\ge 0}$ of non-negative functions on $\R_+$ 
determined by the equation 
\be 
\frac{\partial \psi}{\partial t}(t,\lambda)=R(\psi(t,\lambda)), 
\qquad \psi(0,\lambda)=\lambda 
                                             \label{2.4} 
\ee 
with $\lambda\ge 0$ being arbitrary. 
Let $T_t$ be the semigroup of the CBI-process, 
and for every $\lambda\ge 0$ define a function 
$f_{\lambda}$ on $\R_+$ by 
$f_{\lambda}(x)=e^{-\lambda x}$. 
Then by Theorem 1.1 in \cite{KW71}  
\be 
T_tf_{\lambda}(x)=\exp\left(-x\psi(t,\lambda)
-\int_0^tF(\psi(s,\lambda))ds\right), 
\qquad t, \lambda\ge 0. 
                                             \label{2.5} 
\ee 
Ogura's formula (Eq.(1.12) with $\alpha=0$ in \cite{Og70}) 
for a unique stationary distribution, say $\nu$, of this process is 
\be 
\int_{\R+}f_{\lambda}(x)\nu(dx)= 
\exp\left(-\Phi(\lambda)\right), 
\qquad \lambda\ge 0,   
                                             \label{2.6} 
\ee 
provided the `Laplace exponent' $\Phi$ given by 
\be 
\Phi(\lambda)
=-\int_0^{\lambda}\frac{F(u)}{R(u)}du 
                                             \label{2.7} 
\ee 
is finite for all $\lambda>0$. Conversely, 
if the CBI-process has a stationary distribution, 
then $\Phi(\lambda)<\infty$ for all $\lambda>0$ 
and (\ref{2.6}) holds. (See Lemma 2.1 below for the proof.) 

As shown in \cite{Og70}, the constant $b$ plays an important role 
in studying ergodic properties of the process. 
For instance, under the assumptions $b>0$ and that 
both $R$ and $F$ are analytic at $\lambda=0$, 
the spectral representation of Theorem 3.1 in \cite{Og70} implies 
in particular that $0,b,2b,\ldots$ form 
the discrete spectrum of $-L$. 
Assuming the finiteness of $\Phi$ only, we will see below 
the convergence of the transition function as $t\to\infty$. 
In such a case, the stationary distribution 
$\nu$ is necessarily infinitely divisible 
for the reason mentioned in Introduction, 
and therefore $\Phi$ is expressed uniquely in the form 
\be 
\Phi(\lambda) 
=q\lambda +\int_0^{\infty}(1-e^{-\lambda y})\Lambda(dy) 
                                             \label{2.8} 
\ee 
for some $q\ge 0$ (the `translation term') and 
measure $\Lambda$ (called the L\'evy measure) on $(0,\infty)$ 
such that $\int_{0}^{\infty}\min\{1,y\}\Lambda(dy)<\infty$. 
(See e.g. \S 51 of \cite{Sato}.) 
Obviously $q$ is interpreted as the infimum of the support of $\nu$. 
The condition that $\Phi(1)<\infty$ is sufficient to guarantee  
that $\Phi(\lambda)<\infty$ for every $\lambda>0$ 
since by integration by parts 
\be 
R(\lambda)= -a\lambda^2-b\lambda
-\lambda\int_0^{\infty}\left(1-e^{-\lambda y}\right)
\wt{n_1}(dy)<0,  \qquad \lambda>0, 
                                             \label{2.9} 
\ee 
where $\wt{n_1}(dy)=n_1([y,\infty))dy$. 
Incidentally, we remark that 
$(0,\infty)\ni\lambda\mapsto -R(\lambda)/\lambda$ defines 
a Bernstein function with characteristic triplet 
$(b,a,\wt{n_1})$ in the terminology of \cite{SSV} 
(Chap.3, Theorem 3.2). Let $P_t(x,dy)$ denote 
the transition function of the CBI-process. 
The following lemma gives basic observations concerning ergodicity  
and can be deduced from the results announced in \cite{Pinsky}. 
The proof was given in \cite{Li}. 
(See Theorem 3.20 and Corollary 3.21 there.) 
\begin{lm}
(i) Assume that $\Phi(1)<\infty$ and 
let $\nu$ satisfy (\ref{2.6}). 
Then, for each $x\in\R_+$,  
$P_t(x,\cdot)\to \nu$ weakly as $t\to\infty$.  \\ 
(ii) Suppose that the CBI-process has a stationary distribution, 
then $\Phi(1)<\infty$. \\ 
(iii) If $b>0$, then $\Phi(1)<\infty$. 
\end{lm}
In our discussion, a CBI-process is said to be ergodic 
if it has a (unique) stationary distribution, 
or equivalently $\Phi(1)<\infty$. 
The next proposition concerns not only the translation term $q$ 
of the stationary distribution in the ergodic case 
but also the infimum, denoted by $q(t,x)$, 
of the support of $P_t(x,\cdot)$. 
Put $c=\int_0^{\infty}\wt{n_1}(dy)=\int_0^{\infty}yn_1(dy)$. 
\begin{pr}
(i) If $a>0$ or $c=\infty$, 
then $q(t,x)=0$ for any $x\in\R_+$ and $t>0$. 
Under the additional condition that $\Phi(1)<\infty$, 
it holds that $q=0$. 
\\ 
(ii) If $a=0$ and $0<b+c<\infty$, then for any $x\in\R_+$ and $t>0$ 
\be  
q(t,x)= x e^{-t(b+c)} + \frac{\delta}{b+c}\left(1-e^{-t(b+c)}\right). 
                                             \label{2.11} 
\ee  
Suppose, in addition, that $\Phi(1)<\infty$. Then $q=\delta/(b+c)$. 
\end{pr}
{\it Proof.}~ 
As for $q(t,x)$, the calculation is based on (\ref{2.5}) combined with 
a general fact that 
the infimum of the support of a probability measure $\mu$ on $\R_+$ 
is identified with the `low temperature limit' 
$-\lim_{\lambda\to\infty}(d/{d\lambda})\log\int f_{\lambda}d\mu$. 
Thus 
\be 
q(t,x)= \lim_{\lambda\to\infty}
\left[x\frac{\partial \psi}{\partial \lambda}(t,\lambda)
+\int_0^tF'(\psi(s,\lambda))
\frac{\partial \psi}{\partial \lambda}(s,\lambda)ds\right]. 
                                             \label{2.12} 
\ee 
We introduce an auxiliary function 
$R_0(\lambda)=-R(\lambda)/\lambda$, 
which is positive and increasing for any $\lambda>0$. 
It follows from (\ref{2.4}) that for any $t>0$ and $\lambda>0$  
\be 
t=-\int_{\psi(t,\lambda)}^{\lambda}\frac{du}{R(u)}
=\int_{\psi(t,\lambda)}^{\lambda}\frac{du}{uR_0(u)}.  
                                             \label{2.13} 
\ee 
By differentiating this identity in $\lambda$ 
\be 
\frac{\partial \psi}{\partial \lambda}(t,\lambda) 
= \frac{R(\psi(t,\lambda))}{R(\lambda)} 
= \frac{\psi(t,\lambda)}{\lambda}\cdot 
\frac{R_0(\psi(t,\lambda))}{R_0(\lambda)} 
\in (0,1] 
                                             \label{2.14} 
\ee 
since $\psi(t,\lambda)\le \lambda$. 
Also, noting that $R_0(u)\in[R_0(\psi(t,\lambda)),R_0(\lambda)]$ 
for any $u\in [\psi(t,\lambda), \lambda]$, 
one can deduce from (\ref{2.13}) 
\be 
e^{-tR_0(\lambda)}
\le \frac{\psi(t,\lambda)}{\lambda} 
\le e^{-tR_0(\psi(t,\lambda))}. 
                                             \label{2.15} 
\ee 
(i) To prove that $q(t,x)=0$, we only need to show that 
$\lim_{\lambda \to \infty}\frac{\partial \psi}{\partial \lambda}(t,\lambda)=0$ 
for each $t>0$. Indeed, noting that $F'$ is decreasing 
and that $\psi(s,\lambda)$ is decreasing in $s$ and 
increasing in $\lambda$ by (\ref{2.14}), 
we have for any $\lambda>1$ 
\[ 
0\le 
\int_0^tF'(\psi(s,\lambda))
\frac{\partial \psi}{\partial \lambda}(s,\lambda)ds 
\le F'(\psi(t,1))\int_0^t
\frac{\partial \psi}{\partial \lambda}(s,\lambda)ds.  
\] 
By combining (\ref{2.14}) with (\ref{2.15}) 
\[ 
0\le 
\frac{\partial \psi}{\partial \lambda}(t,\lambda) 
\le 
e^{-tR_0(\psi(t,\lambda))}
\frac{R_0(\psi(t,\lambda))}{R_0(\lambda)}
\le \frac{1}{tR_0(\lambda)}, 
\] 
which tends to 0 as $\lambda \to \infty$ 
since the assumption implies that $R_0(\lambda)\to \infty$. 
Consequently $\frac{\partial \psi}{\partial \lambda}(t,\lambda)$ 
converges to 0 boundedly and 
by virtue of (\ref{2.12}) $q(t,x)=0$. 

In the case where $\Phi(1)<\infty$, 
$q=\lim_{\lambda\to\infty}\Phi'(\lambda)
=\lim_{\lambda\to\infty}F(\lambda)/(\lambda R_0(\lambda))$. 
It is easy to see that $\lim_{\lambda\to\infty}F(\lambda)/\lambda=\delta$. 
So we conclude that $q=0$.  \\ 
(ii) It is obvious that the proof of (\ref{2.11}) 
can be reduced to showing the following 
two asymptotics; as $\lambda\to\infty$  
\[ 
\frac{\partial \psi}{\partial \lambda}(s,\lambda) 
\to e^{-(b+c)s} 
\quad \mbox{for each} \ s>0 
\] 
and 
\[ 
F'(\psi(s,\lambda)) \to \delta 
\quad \mbox{locally boundedly in} \ s\ge 0. 
\] 
Observe that $R_0(\lambda)\to b+c\in(0,\infty)$ by the assumption. 
Hence the first inequality in (\ref{2.15}) implies that 
$\psi(s,\lambda)\to\infty$. By (\ref{2.15}) again 
we have $\psi(s,\lambda)/\lambda \to \exp(-(b+c)s)$ 
and thus (\ref{2.14}) proves the first asymptotics. 
The second one is a consequence of the following estimate; 
for any $s\in[0,t]$ 
\[ 
|F'(\psi(s,\lambda))-\delta|
= \int_0^{\infty}y e^{-\psi(s,\lambda)y}n_2(dy)
\le \int_0^{\infty}y e^{-\psi(t,\lambda)y}n_2(dy). 
\] 
Therefore (\ref{2.11}) has been established. 

The last part of the assertion (ii) can be shown 
by $F(\lambda)/\lambda\to \delta$ and 
$R_0(\lambda)\to b+c$ together. 
The proof of the proposition is complete. 
\qed 
\bigskip 

It is worth mentioning the diffusion case, namely the case 
where $a,b>0$ and $n_1\equiv 0 \equiv n_2$. 
Because of (\ref{2.3}) we have $\delta>0$, and by (\ref{2.7}) 
the corresponding CBI-process, known also as the CIR model, 
has a unique stationary distribution with Laplace exponent 
\[ 
\Phi(\lambda)
=\frac{\delta}{a}\log\left(1+\frac{a}{b}\lambda\right) 
=\frac{\delta}{a}\int_0^{\infty}(1-e^{-\lambda y})\frac{e^{-by/a}}{y}dy. 
\] 
It is a gamma distribution 
with parameter $(\delta/a,b/a)$, which has, 
by definition, the density proportional to 
$x^{\delta/a-1}\exp(-bx/a)$. 
This stationary distribution is reversible. 
In other words, the associated bilinear form (\ref{1.1}) is symmetric. 

It would be natural to ask if there is any other case 
which admits a nondegenerate reversible distribution. 
The following theorem, the main result of this section, 
gives a negative answer to this question. 
\begin{th}
If the CBI-process with generator (\ref{2.1}) 
has a nondegenerate reversible distribution, 
then the process coincides with a CIR model. 
\end{th}
Before proving this theorem we prepare a simple lemma. 
In what follows, the notation $\lg f\rg$ 
or $\lg f(x)\rg$ will stand for 
the integral $\int_{\R_+}f(x)\nu(dx)$ 
with respect to the stationary distribution $\nu$ 
of an ergodic CBI-process. 
\begin{lm}
Let $\lambda, \mu\ge 0$ be arbitrary. \\ 
(i) For any $x\in \R_+$ 
\be 
-Lf_{\lambda}(x)= (R(\lambda)x+F(\lambda))f_{\lambda}(x).  
                                             \label{2.16} 
\ee 
(ii) If the CBI-process has a stationary distribution, then 
\be 
\lg (-L)f_{\lambda}\cdot f_{\mu}\rg
= R(\lambda)(\Phi'(\lambda+\mu)-\Phi'(\lambda))
\lg f_{\lambda+\mu}\rg. 
                                             \label{2.17} 
\ee 
\end{lm}
{\it Proof.}~ 
(\ref{2.16}) is verified by direct calculations. 
Using it, we have 
\be  
\lg (-L)f_{\lambda}\cdot f_{\mu}\rg
=R(\lambda)\lg xe^{-(\lambda+\mu)x}\rg 
+ F(\lambda) \lg e^{-(\lambda+\mu)x}\rg. 
                                             \label{2.18} 
\ee  
Also, (\ref{2.6}) and (\ref{2.7}) give 
$\lg xe^{-(\lambda+\mu)x}\rg
=\Phi'(\lambda+\mu)\lg e^{-(\lambda+\mu)x}\rg$ 
and $F(\lambda)=-R(\lambda)\Phi'(\lambda)$, respectively. 
(\ref{2.17}) follows by plugging these equalities into (\ref{2.18}). 
\qed 

\bigskip 

\noindent 
{\it Proof of Theorem 2.3.}~
By the assumption, we have 
the symmetry of the Dirichlet form. In particular, 
$\lg(-L)f_{\lambda}\cdot f_{\mu}\rg=\lg(-L)f_{\mu}\cdot f_{\lambda}\rg$  
for all $\lambda, \mu>0$. By virtue of Lemma 2.4 (ii), this becomes 
\[ 
R(\lambda)(\Phi'(\lambda+\mu)-\Phi'(\lambda))
=R(\mu)(\Phi'(\lambda+\mu)-\Phi'(\mu)). 
\] 
Because of (\ref{2.7}), the above equality is rewritten into 
\be 
(R(\lambda)-R(\mu))F(\lambda+\mu)
=(F(\lambda)-F(\mu))R(\lambda+\mu). 
                                             \label{2.19} 
\ee 
Differentiating in $\lambda$ yields 
\[ 
R'(\lambda)F(\lambda+\mu)+(R(\lambda)-R(\mu))F'(\lambda+\mu)
=F'(\lambda)R(\lambda+\mu)+(F(\lambda)-F(\mu))R'(\lambda+\mu). 
\] 
By interchanging the roles of $\lambda$ and $\mu$ 
\[ 
R'(\mu)F(\lambda+\mu)+(R(\mu)-R(\lambda))F'(\lambda+\mu)
=F'(\mu)R(\lambda+\mu)+(F(\mu)-F(\lambda))R'(\lambda+\mu). 
\] 
Summing up the above two equalities, we arrive at 
\be 
(R'(\lambda)+R'(\mu))F(\lambda+\mu)
=(F'(\lambda)+F(\mu))R(\lambda+\mu). 
                                             \label{2.20} 
\ee 
Now let $\lambda\ne\mu$. Then 
(\ref{2.19}) and (\ref{2.20}) together with 
Cauchy's mean value theorem imply further that 
for some $\xi$ between $\lambda$ and $\mu$ 
\be 
\frac{R'(\lambda)+R'(\mu)}{F'(\lambda)+F'(\mu)}
=\frac{R(\lambda)-R(\mu)}{F(\lambda)-F(\mu)}
=\frac{R'(\xi)}{F'(\xi)}. 
                                             \label{2.21} 
\ee 
Here it should be noted that by (\ref{2.3})  
\be 
F'(u)=\delta+\int_0^{\infty}ye^{-uy}n_2(dy)>0 
                                             \label{2.22} 
\ee 
is convex and 
\[ 
R'(u)=-2au-b-\int_0^{\infty}y(1-e^{-uy})n_1(dy)<0
\] 
is convex. 
Therefore, (\ref{2.21}) is possible only in the case where 
neither $F'$ nor $-R'$ is strictly convex. 
Consequently, both $n_1$ and $n_2$ must vanish. 
So (\ref{2.22}) shows that $\delta>0$, and 
the positivity of $b$ is necessary for $\Phi(\lambda)$ to be finite. 
Moreover, $a$ must be positive also since 
otherwise $\Phi(\lambda)=\delta\lambda/b$ implying 
that the stationary distribution is concentrated at $\delta/b$. 
The proof of Theorem 2.3 is complete. 
\qed 

\bigskip 

Here are concrete examples of CBI-processes 
with non-reversible stationary distributions. 

\medskip 

\noindent 
{\it Example 2.1}~
(i) This example is taken from Example 4.2 of \cite{Og70}. 
Given $0<\alpha<\beta<1$, set $a=b=\delta=0$, 
\[ 
n_1(dy)=\frac{\alpha(\alpha+1)}{\Gamma(1-\alpha)}
\cdot\frac{dy}{y^{2+\alpha}} 
\quad \mbox{and} \quad 
n_2(dy)=\frac{\beta}{\Gamma(1-\beta)}
\cdot\frac{dy}{y^{1+\beta}}. 
\] 
Then $R(\lambda)=-\lambda^{1+\alpha}$ 
and $F(\lambda)=\lambda^{\beta}$. Therefore, 
$\Phi(\lambda)=\lambda^{\beta-\alpha}/(\beta-\alpha)$, 
the Laplace exponent of a $(\beta-\alpha)$-stable distribution on $\R_+$. 
(\ref{2.17}) gives 
\[ 
\lg (-L)f_{\lambda}\cdot f_{\mu}\rg
= 
\lambda^{\beta}
\left\{1-\left(\frac{\lambda}{\lambda+\mu}\right)^{1-(\beta-\alpha)}\right\}
e^{-\Phi(\lambda+\mu)}. 
\] 
(ii) 
Given $0<\alpha<1$ and $\kappa\ge 0$, define $a=0$, 
$b=\kappa^{\alpha}$, $\delta=1$, $n_2\equiv 0$ and 
\be 
n_1(dy)=-\frac{\alpha}{\Gamma(1-\alpha)}
\left(\frac{e^{-\kappa y}}{y^{1+\alpha}}\right)'dy 
=\frac{\alpha}{\Gamma(1-\alpha)}(\kappa y+1+\alpha)
\frac{e^{-\kappa y}}{y^{2+\alpha}}dy. 
                                             \label{2.23} 
\ee 
With these choices 
$R(\lambda)=-\lambda(\lambda+\kappa)^{\alpha}$ 
and $F(\lambda)=\lambda$, which together lead to 
\be 
\Phi(\lambda)=
\frac{1}{1-\alpha}
\left[(\lambda+\kappa)^{1-\alpha}-\kappa^{1-\alpha}\right] 
=\frac{1}{\Gamma(\alpha)}
\int_0^{\infty}(1-e^{-\lambda y}) 
\frac{e^{-\kappa y}}{y^{1+(1-\alpha)}}dy. 
                                             \label{2.24} 
\ee 
See e.g. \cite{Hougaard} for information of 
the corresponding distribution. 
Note that for $\kappa=0$ the stationary distribution is a 
$(1-\alpha)$-stable distribution on $\R_+$, and that 
as $\alpha\uparrow 1$, 
$\Phi(\lambda)$ tends to $\log(1+\lambda/\kappa)$, 
the Laplace exponent of a gamma distribution, 
provided that $\kappa>0$. By (\ref{2.17}) we have 
\[ 
\lg (-L)f_{\lambda}\cdot f_{\mu}\rg
= 
\lambda
\left\{1-\left(\frac{\lambda+\kappa}
{\lambda+\mu+\kappa}\right)^{\alpha}\right\}
e^{-\Phi(\lambda+\mu)}. 
\] 

\medskip 

Since the class of CBI-processes studied so far 
seems too wide for one to obtain further consequences 
which are useful for our purpose, 
we will be obliged to make an additional restriction 
in the subsequent sections. 
In this regard, it must be remarked that 
the condition that $n_1\equiv 0$ makes 
the correspondence between $(a,b,n_2,\delta)$ and 
$(q,\Lambda)$ in (\ref{2.8}) 
too simple as will be seen from the general observation below. 
\begin{lm}
Let $a,b\ge 0$ and 
suppose that a measure $n_2$ on $(0,\infty)$ is non-zero. 
For $\lambda>0$, set 
\[ 
\Phi(\lambda)
=\int_0^{\lambda}
\frac{\ds{\int_0^{\infty}(1-e^{-uy})n_2(dy)}}{au^2+bu}du
=\int_0^{\lambda}
\frac{\ds{\int_0^{\infty}e^{-uy}\wt{n_2}(dy)}}{au+b}du. 
\] 
Then for each $\lambda>0$ 
\[  
\Phi(\lambda)=
\left\{
\begin{array}{ll} 
\infty 
& (b=0) \\ 
\ds{\int_0^{\infty}(1-e^{-\lambda y})\frac{\wt{n_2}(dy)}{by}} 
& (a=0, \ b>0) \\ 
\ds{\int_0^{\infty}(1-e^{-\lambda y})
\left(\int_0^ye^{bz/a}\wt{n_2}(dz)\right)\frac{e^{-by/a}}{ay}dy} 
& (a>0, \ b>0). 
\end{array}
\right. 
\] 
\end{lm}
The proof requires only `Fubini calculus'  
and so is left to the reader. 
In the light of this lemma, 
we shall proceed under the additional hypothesis that $n_2\equiv 0$.

\section{Estimating the sector constant}
\setcounter{equation}{0} 

The main subject of this section is 
the estimation of the Dirichlet form. 
As announced in Introduction, we now show in a general setting that 
${\rm Sect}(\cE)>1$ holds for any non-symmetric Dirichlet form $\cE$. 
(An explicit lower bound for ${\rm Sect}(\cE)$ 
will be discussed at the end of this section.) 
\begin{pr}
Suppose that the bilinear form $\cE$ in (\ref{1.1}) is 
associated with a conservative Markov process 
with generator $L$ and a stationary distribution $\nu$. 
If $\cE$ is non-symmetric, then ${\rm Sect}(\cE)>1$. 
\end{pr}
{\it Proof.}~ 
We may assume that ${\rm Sect}(\cE)<\infty$. 
Equivalently, suppose that (\ref{1.2}) holds for some $C<\infty$. 
By non-symmetry there exist $f,g\in D(L)$ such that 
$\cE(f,g)>\cE(g,f)$. 
This implies $\cE(f,f)>0$ since otherwise 
(\ref{1.2}) leads to the contradiction 
that $\cE(f,g)=\cE(g,f)(=0)$. 
It is straightforward to see that 
\be 
\lim_{t\to 0}
\frac{1}{t}
\left(
\frac{\cE(f,f+t g)^2}
{\cE(f,f)\cE(f+t g,f+t g)}-1\right)
=\frac{\cE(f,g)-\cE(g,f)}{\cE(f,f)}>0, 
                                             \label{3.1} 
\ee 
and hence $\cE(f,f+t g)^2> \cE(f,f)\cE(f+t g,f+t g)$ 
for $t>0$ small enough. This shows that ${\rm Sect}(\cE)>1$. 
\qed 

\medskip 

We now turn to discussing the bilinear forms $\cE$ 
associated with ergodic CBI-processes. 
The symmetric part $\wt{\cE}(f,g):=(\cE(f,g)+\cE(g,f))/2$  
has an expression of the form 
\begin{eqnarray}   
\wt{\cE}(f,g)
& = & 
a\lg xf'(x)g'(x)\rg 
+\frac{1}{2}\lg x\int n_1(dy)(f(x+y)-f(x))(g(x+y)-g(x))\rg  
\nonumber \\ 
& & +\frac{1}{2}\lg \int n_2(dy)(f(x+y)-f(x))(g(x+y)-g(x))\rg. 
                                             \label{3.2} 
\end{eqnarray} 
Here and in what follows, the domain of integration may be suppressed 
as long as it is $(0,\infty)$. 
(\ref{3.2}) can be verified by calculating 
$\Gamma(f,g):=(L(fg)-Lf\cdot g-f\cdot Lg)/2$ 
since $\wt{\cE}(f,g)=\lg \Gamma(f,g) \rg$. 
The main task in the remainder of this section is to give 
an upper bound of ${\rm Sect}(\cE)$ for a class of 
non-reversible CBI-processes. 

In the rest of the paper, we make the restriction that $n_2\equiv 0$ 
and call such processes 
continuous state branching processes  
with continuous immigration, 
henceforth abbreviated as CBCI-processes. 
Thus, $F(\lambda)=\delta\lambda$ with some $\delta>0$. 
This condition seems crucial in the subsequent argument, 
in particular, in showing an integration by parts formula 
described in Proposition 3.2 below. 
The notation $n$ is used instead of $n_1$ and 
thus a measure $n$ on $(0,\infty)$ is assumed to satisfy 
$\int \min\{y^2,y\}n(dy)<\infty$ according to (\ref{2.2}). 
In the discussion below, we shall suppose 
the existence of a unique stationary distribution 
and introduce a one-parameter family of the convolution semigroup 
$\{\nu_{\delta}:=\nu^{\ast\delta}\}_{\delta>0}$ 
with $\nu$ having the Laplace exponent 
\be 
\Phi(\lambda)
=\int_0^{\lambda}\frac{du}{\ds{au+b+\int(1-e^{-uy})\wt{n}(dy)}}
=q\lambda+\int(1-e^{-\lambda y})\Lambda(dy), 
                                            \label{3.3} 
\ee 
where $q\ge 0$ and $\Lambda$ is a L\'evy measure on $(0,\infty)$. 
Accordingly, $\delta\Phi$ is the Laplace exponent of $\nu_{\delta}$, 
which is a unique stationary distribution of the process with generator 
\begin{eqnarray*}
L_{\delta}f(x) 
&:=& axf''(x)-bxf'(x) 
+x\int\left[f(x+y)-f(x)-yf'(x)\right]n(dy) \nonumber  \\ 
& & +\delta f'(x). 
\end{eqnarray*}
We call it the CBCI-process with quadruplet $(a,b,n,\delta)$. 
This subclass of CBI-processes is 
one-dimensional version of the model 
discussed in \cite{Stannat03} and \cite{Stannat05}. 
We emphasize that an explicit formula 
for the L\'evy density $d\Lambda/dy$ 
was obtained in Lemma 2.5 of \cite{Stannat05} 
under the additional hypothesis 
that $a,b>0$ and $c=\int \wt{n}(dy)<\infty$. 
(In this case $q=0$ by Proposition 2.2 (i).) 
An analogue of that formula is available also 
in the case where $a=0$ and $0<b+c<\infty$. 
Indeed, by differentiating (\ref{3.3}) 
\[ 
\Phi'(\lambda) 
= \frac{1}{\ds{b+c-\int e^{-\lambda y}\wt{n}(dy)}} 
= \frac{1}{b+c} 
+\sum_{N=1}^{\infty}\frac{1}{b+c}
\left(\frac{1}{b+c}\int e^{-\lambda y}\wt{n}(dy)\right)^N,  
\] 
which shows that $q=1/(b+c)$ and 
\be 
\Lambda(dy)=\frac{1}{y}\sum_{N=1}^{\infty}\frac{1}{(b+c)^{N+1}} 
\wt{n}^{*N}(dy). 
                                             \label{3.4}
\ee 
The notation $\lg \cdot\rg_{\delta}$ will stand for 
the integral with respect to $\nu_{\delta}$ 
and the associated bilinear form is denoted 
by $\cE^{\delta}$, namely 
$\cE^{\delta}(f,g)=\lg (-L_{\delta})f \cdot g \rg_{\delta}$. 
Let $\cF_0$ be the linear hull of $\{f_{\lambda}: \lambda\ge 0\}$. 
As remarked after Theorem 1.1$'$ in \cite{KW71}, 
$\cF_0$ is a core of the generator of the CBI-process. 
\begin{pr}
Let $\Phi$ be defined by the first equality in (\ref{3.3}) and 
suppose that $\Phi(1)<\infty$. Then for each $\delta>0$ 
\begin{eqnarray}   
\cE^{\delta}(f,g)  
& = & a\lg xf'(x)g'(x)\rg_{\delta}                 \nonumber   \\ 
&  &+\lg x\int \wt{n}(dy)f'(x+y)(g(x+y)-g(x))\rg_{\delta}, 
\quad    f,g\in\cF_0. 
                                             \label{3.5} 
\end{eqnarray}   
\end{pr}
{\it Proof.}~ 
It suffices to show (\ref{3.5}) for 
$f=f_{\lambda}$, $g=f_{\mu}$ with $\lambda,\mu>0$ being arbitrary.  
We begin with a version of (\ref{2.17}) in Lemma 2.4: 
\be 
\lg (-L_{\delta})f_{\lambda} \cdot f_{\mu} \rg_{\delta} 
= \delta\Phi'(\lambda+\mu)\lg f_{\lambda+\mu}\rg_{\delta} 
\left(R(\lambda)+\frac{\lambda}{\Phi'(\lambda+\mu)}\right),   
                                             \label{3.6} 
\ee 
where 
$R(\lambda)= -a\lambda^2-b\lambda
-\lambda\int \wt{n}(dy)(1-e^{-\lambda y})$. 
Observing from (\ref{3.3}) that 
\[ 
\frac{1}{\Phi'(\lambda+\mu)}
=a(\lambda+\mu)+b+\int\wt{n}(dy)(1-e^{-(\lambda+\mu)y}),  
\] 
we have 
\be  
R(\lambda)+\frac{\lambda}{\Phi'(\lambda+\mu)}
=a\lambda\mu+\lambda\int\wt{n}(dy)
\left(e^{-\lambda y}-e^{-(\lambda+\mu)y}\right). 
                                             \label{3.7} 
\ee 
Since $\delta\Phi'(\lambda+\mu)\lg f_{\lambda+\mu}\rg_{\delta}
=\lg xf_{\lambda+\mu}(x)\rg_{\delta}$, 
we get by plugging (\ref{3.7}) into (\ref{3.6}) 
\begin{eqnarray} 
\lg (-L_{\delta})f_{\lambda} \cdot f_{\mu} \rg_{\delta} 
& = & \lg xf_{\lambda+\mu}(x)\rg_{\delta}
\left(a\lambda\mu+\lambda\int\wt{n}(dy)
\left(e^{-\lambda y}-e^{-(\lambda+\mu)y}\right)\right) 
                                             \label{3.8}  \\ 
& = & 
a\lg x{f_{\lambda}}'(x){f_{\mu}}'(x)\rg_{\delta} 
+ \lg x\int\wt{n}(dy)(-\lambda)e^{-\lambda(x+y)}
\left(e^{-\mu(x+y)}-e^{-\mu x}\right)\rg_{\delta}. 
\nonumber 
\end{eqnarray} 
This coincides with the right  side of (\ref{3.5}) 
with $f_{\lambda}$ and $f_{\mu}$ in place of $f$ and $g$, respectively. 
\qed 

\bigskip 

The integration by parts formula (\ref{3.5}) would be 
interesting in its own right and applicable in other contexts. 
We here use it for the purpose of the sector constant estimate. 
While it is obvious from (\ref{3.2}) that the first term on 
the right side of (\ref{3.5}) is dominated by 
${\cE^{\delta}(f,f)}^{1/2}{\cE^{\delta}(g,g)}^{1/2}$, 
the main difficulty in handling the second term comes from 
the fact that we have few information on 
the distribution function (or the density function) 
of the stationary distribution $\nu_{\delta}$. 
We overcome this by a strategy similar to 
that taken in \cite{Stannat05} for the proof of 
a Poincar\'e type inequality. That is, 
we show first that an estimate 
we want for the second term to satisfy 
can reduce to an analogous one for 
the L\'evy measure $\Lambda$, 
and then give a sufficient condition 
for the reduced estimate to hold. 
So we shall be concerned with the bilinear forms 
$B_{\delta}(\delta>0)$ and $B_{0}$ 
on $\cF_0\times \cF_0$ defined by 
\[ 
B_{\delta}(f,g)
=\lg x\int\wt{n}(dy)f'(x+y)(g(x+y)-g(x))\rg_{\delta}, 
\qquad \delta>0 
\] 
and 
\[ 
B_{0}(f,g)
=\int \Lambda(dx)x\int\wt{n}(dy)f'(x+y)(g(x+y)-g(x)),  
\] 
respectively. The first step will be done in the next theorem, 
the main technical result of this section. 
\begin{th}
Let $\Phi$ be defined by the first equality in (\ref{3.3}). 
Suppose that $\Phi(1)<\infty$ and that 
a L\'evy measure $\Lambda$ on $(0,\infty)$ 
satisfies (\ref{3.3}) with $q=0$. 
Then, for any fixed $0<C<\infty$, 
the following two conditions are equivalent to each other: \\ 
(i) For all $\delta>0$ 
\be 
B_{\delta}(f,g)^2
\le C\lg xf'(x)^2\rg_{\delta} 
\lg x\int n(dy)(g(x+y)-g(x))^2\rg_{\delta}, 
\quad f,g\in\cF_0. 
                                             \label{3.9}
\ee 
(ii) For all $f,g\in\cF_0$ 
\be 
B_{0}(f,g)^2
\le C\int \Lambda(dx) xf'(x)^2\cdot 
\int \Lambda(dx) x\int n(dy)(g(x+y)-g(x))^2.  
                                             \label{3.10} 
\ee  
If, in addition, $a>0$ and (\ref{3.9}) holds for some $\delta>0$,  
then ${\rm Sect}(\cE^{\delta})\le 1+\sqrt{2C/a}$. 
\end{th}
{\it Proof.}~ 
The implication (i)$\Longrightarrow$(ii) follows 
immediately by observing that as $\delta\downarrow 0$ 
\be 
\delta^{-1}\lg x h_i(x)\rg_{\delta} 
= \int z\Lambda(dz)\lg h_i(x+z)\rg_{\delta}
\to \int z\Lambda(dz) h_i(z) 
                                             \label{3.11} 
\ee 
for each $i\in\{1,2,3\}$, 
where 
$h_1(x)=\int\wt{n}(dy)f'(x+y)(g(x+y)-g(x))$, 
$h_2(x)=f'(x)^2$ and 
$h_3(x)=\int n(dy)(g(x+y)-g(x))^2$ with 
$f,g\in\cF_0$ being arbitrary. 
In (\ref{3.11}) we have applied 
the Palm formula for the underlying Poisson random measure 
to get the equality, and then used 
the fact that $\nu_{\delta}$ tends weakly 
to the delta distribution at $0$. 
(Alternatively, the equality can be verified 
directly for $h_i=f_{\lambda}$ with $\lambda>0$ and extended easily. 
See Lemma 3.2 in \cite{Stannat05}.) 

Next, assume that (\ref{3.10}) holds for all $f,g\in\cF_0$. 
We must show (\ref{3.9}) for every $\delta>0$. 
Let $f$ and $g$ be given as finite sums of the form 
$f=\sum_i c_i f_{\lambda_i}$ and $g=\sum_j d_j f_{\mu_j}$, 
respectively, 
where $c_i,d_j\in\R, \lambda_i,\mu_j\ge 0$. 
In view of (\ref{3.8}) 
\begin{eqnarray*} 
B_{\delta}(f_{\lambda},f_{\mu}) 
& = & 
\lg xf_{\lambda+\mu}(x)\rg_{\delta}
\lambda \int\wt{n}(dy)
\left(e^{-\lambda y}-e^{-(\lambda+\mu)y}\right) \\ 
& = & 
\delta e^{-\delta \Phi(\lambda+\mu)} 
\Phi'(\lambda+\mu)\int\wt{n}(dy)\lambda 
\left(e^{-\lambda y}-e^{-(\lambda+\mu)y}\right)  
\end{eqnarray*} 
for any $\lambda,\mu\ge 0.$ 
Therefore, by bilinearity $B_{\delta}(f,g)$ equals 
\begin{eqnarray} 
\lefteqn{\delta \sum_{i,j}c_id_j 
e^{-\delta \Phi(\lambda_i+\mu_j)} 
\Phi'(\lambda_i+\mu_j)\int\wt{n}(dy)\lambda_i 
\left(e^{-\lambda_i y}-e^{-(\lambda_i+\mu_j)y}\right)} \nonumber \\ 
& = & 
\delta \sum_{i,j}c'_id'_j 
e^{\delta \Phi(\lambda_i,\mu_j)} 
\int\Lambda(dx)xe^{-(\lambda_i+\mu_j)x}
\int\wt{n}(dy)\lambda_i 
\left(e^{-\lambda_i y}-e^{-(\lambda_i+\mu_j)y}\right), 
                                             \label{3.12} 
\end{eqnarray} 
where $c'_i=c_i\exp(-\delta\Phi(\lambda_i))$, 
$d'_j=d_j\exp(-\delta\Phi(\mu_j))$ and 
\[ 
\Phi(\lambda,\mu)
=\Phi(\lambda)+\Phi(\mu)-\Phi(\lambda+\mu)
=\int\Lambda(dz)(1-e^{-\lambda z})(1-e^{-\mu z}), 
\qquad \lambda,\mu\ge 0. 
\] 
Substituting the expansion 
\[ 
e^{\delta \Phi(\lambda_i,\mu_j)} 
=\sum_{N=0}^{\infty}\frac{\delta^N}{N!}
\int\Lambda(dx_1)\cdots\int\Lambda(dx_N)
\prod_{k=1}^N(1-e^{-\lambda_ix_k})
\prod_{l=1}^N(1-e^{-\mu_jx_l}) 
\] 
into (\ref{3.12}) leads to 
\begin{eqnarray*}
B_{\delta}(f,g) 
& =  & 
\delta \sum_{N=0}^{\infty}\frac{\delta^N}{N!}
\int\Lambda(dx_1)\cdots\int\Lambda(dx_N)\int\Lambda(dx)x \\ 
& & \int\wt{n}(dy){f_{x_1,\ldots,x_N}}'(x+y)
\left(g_{x_1,\ldots,x_N}(x+y)
-g_{x_1,\ldots,x_N}(x)\right), 
\end{eqnarray*} 
where 
\[ 
f_{x_1,\ldots,x_N}(x)=
\sum_ic'_i\prod_{k=1}^N(1-e^{-\lambda_ix_k})f_{\lambda_i}(x) 
\] 
and 
\[ 
g_{x_1,\ldots,x_N}(x)=
\sum_jd'_j\prod_{l=1}^N(1-e^{-\mu_jx_l})f_{\mu_j}(x) 
\] 
are considered to be elements of $\cF_0$ 
for arbitrarily given $x_1,\ldots,x_N>0$. 
We can apply now (\ref{3.10}) to these functions 
and then use Schwarz's inequality to obtain 
\begin{eqnarray}
B_{\delta}(f,g) 
& =  & 
\delta \sum_{N=0}^{\infty}\frac{\delta^N}{N!}
\int\Lambda(dx_1)\cdots\int\Lambda(dx_N)
B_{0}(f_{x_1,\ldots,x_N},g_{x_1,\ldots,x_N})  
                                             \label{3.13} \\ 
& \le & 
\sqrt{C} \sqrt{Q_{\delta}^{(1)}(f)}\sqrt{Q_{\delta}^{(2)}(g)}, 
                                             \label{3.14} 
\end{eqnarray} 
where 
\[ 
Q_{\delta}^{(1)}(f)
= \delta \sum_{N=0}^{\infty}\frac{\delta^N}{N!}
\int\Lambda(dx_1)\cdots\int\Lambda(dx_N)
\int\Lambda(dx)x {f_{x_1,\ldots,x_N}}'(x)^2 
\] 
and $Q_{\delta}^{(2)}(g)$ is defined to be 
\[ 
\delta \sum_{N=0}^{\infty}\frac{\delta^N}{N!}
\int\Lambda(dx_1)\cdots\int\Lambda(dx_N)
\int\Lambda(dx)x \int n(dy) 
\left(g_{x_1,\ldots,x_N}(x+y)-g_{x_1,\ldots,x_N}(x)\right)^2. 
\] 
But analogous calculations to those 
for the proof of (\ref{3.13}) show that 
\[ 
Q_{\delta}^{(1)}(f)
=\lg xf'(x)^2\rg_{\delta} 
\quad \mbox{and} \quad  
Q_{\delta}^{(2)}(g)
=\lg x\int n(dy)(g(x+y)-g(x))^2\rg_{\delta}.   
\] 
Thus, (\ref{3.14}) proves (\ref{3.9}). 

Lastly, in view of (\ref{3.5}) and (\ref{3.2}) 
with $n_1=n$ and $n_2\equiv 0$, 
the validity of (\ref{3.9}) implies that 
\[ 
\cE^{\delta}(f,g)\le (1+\sqrt{2C/a})
{\cE^{\delta}(f,f)}^{1/2}{\cE^{\delta}(g,g)}^{1/2}, 
\qquad f,g\in \cF_0 
\] 
if $a>0$. This inequality is shown to extend 
to all functions in $D(L_{\delta})$, 
and hence ${\rm Sect}(\cE^{\delta})\le 1+\sqrt{2C/a}$ as desired. 
The proof of Theorem 3.3 is complete. \qed 

\bigskip 

The next step is to seek conditions for (\ref{3.10}) to hold. 
Define 
\[ 
\left\Vert \frac{d\wt{n}}{dn} \right\Vert_{\infty}
=\inf\{r>0:~ 
\wt{n}(dy)\le r n(dy) \ \mbox{in distribution sense} \}
\] 
with convention that $\inf\emptyset =\infty$. 
Clearly, this value is 0 for $n\equiv 0$. 
\begin{cor}
In addition to the assumptions in Theorem 3.3, 
assume that $a>0$ and that there exists a density 
$d \Lambda/dy=:\varphi(y)/y$ such that 
\be 
(\varphi*\wt{n})(z):=\int_0^z 
\varphi(z-y)\wt{n}(dy) \le C_1 \varphi(z), 
\qquad z>0 
                                             \label{3.15} 
\ee 
for some $0\le C_1<\infty$. Then for any $\delta>0$ 
\be 
{\rm Sect}(\cE^{\delta})\le 1
+\sqrt{\frac{2C_1}{a}
\left\Vert \frac{d\wt{n}}{dn} \right\Vert_{\infty}}.  
                                             \label{3.16} 
\ee 
\end{cor}
{\it Proof.}~ 
We may assume that $\Vert d\wt{n}/dn \Vert_{\infty}<\infty$. 
By virtue of Theorem 3.3, 
it is enough to show (\ref{3.10}) with $C=C_1r$ for 
any $r>0$ such that $\wt{n}(dy)\le r n(dy)$. 
Applying Schwarz's inequality, 
we can dominate $B_{0}(f,g)^2$ by    
\begin{eqnarray*}
\lefteqn{
\int \Lambda(dx) x\int \wt{n}(dy)f'(x+y)^2\cdot 
\int \Lambda(dx) x\int \wt{n}(dy)(g(x+y)-g(x))^2}   \\ 
& \le & 
\int dx \varphi(x) \int \wt{n}(dy)f'(x+y)^2 
\cdot r \int \Lambda(dx) x\int n(dy)(g(x+y)-g(x))^2.  
\end{eqnarray*} 
Since by (\ref{3.15}) 
\begin{eqnarray*}
\int dx \varphi(x) \int \wt{n}(dy)f'(x+y)^2 
& = & 
\int dz f'(z)^2(\varphi*\wt{n})(z) \\ 
& \le  & 
C_1 \int dz f'(z)^2\varphi(z)  = 
C_1 \int \Lambda(dz)zf'(z)^2,  
\end{eqnarray*}
the desired inequality is derived. 
\qed 

\bigskip 

A key ingredient to verify (\ref{3.15}) is 
the following fact taken from Eq. (6) 
in the proof of Lemma 2.6 in \cite{Stannat05}. 
(The function $K$ there is identical with $\varphi$ 
in the present paper.) 
\begin{lm}
In addition to the assumptions in Theorem 3.3, 
assume that $a>0$ and $0<\wt{b}:=b+c<\infty$. 
Then the L\'evy measure $\Lambda$ in (\ref{3.3}) has 
a strictly positive density 
$d\Lambda/dy=:\varphi(y)/y$ with $\varphi$ being differentiable. 
Moreover, $\varphi(0):=\lim_{y\downarrow 0}\varphi(y)=1/a$ and  
\be 
(\varphi*\wt{n})(y)=
a\varphi'(y)+\wt{b}\varphi(y), 
\qquad y>0. 
                                             \label{3.17} 
\ee 
\end{lm}
To grasp the validity of (\ref{3.17}), 
it is worth noting that 
taking the Laplace transform of both side of (\ref{3.17}) leads, 
at least at formal level, to the equation equivalent to 
the one derived by differentiating (\ref{3.3}) with $q=0$ 
provided that $\varphi(0)=1/a$. 
Note also that 
$\varphi'(0):=\lim_{y\downarrow 0}\varphi'(y)=-\wt{b}/a^2$ 
is deduced from (\ref{3.17}). 
\begin{pr}
Under the same assumptions and with the same notation 
as in Lemma 3.5, define $V(y)=-\log \varphi(y)$ 
(, so that $\lim_{y\downarrow 0}V'(y)=:V'(0)$ exists). 
If 
\[ 
\sup\{V'(0)-V'(y):~y>0 \}\le C_2 
\] 
for some $0\le C_2<\infty$, then for any $\delta>0$ 
\be  
{\rm Sect}(\cE^{\delta})
\le 1+\sqrt{2C_2\left\Vert \frac{d\wt{n}}{dn} \right\Vert_{\infty}}.   
                                             \label{3.18} 
\ee 
\end{pr}
{\it Proof.}~ 
By letting $y\downarrow 0$ in (\ref{3.17}) 
$a\varphi'(0)+\wt{b}\varphi(0)=0$ or 
$\wt{b}=-a\varphi'(0)/\varphi(0)=aV'(0)$. 
So again by (\ref{3.17}) 
\[ 
(\varphi*\wt{n})(y)
= a\varphi'(y)+a\varphi(y)V'(0)
= a\varphi(y)(-V'(y)+V'(0))
\le aC_2\varphi(y), 
\] 
and thus (\ref{3.15}) with $C_1=aC_2$ holds true. 
Therefore (\ref{3.18}) follows from (\ref{3.16}). 
\qed 

\bigskip 

\noindent 
In the reversible case $n\equiv 0$, 
the function $V$ in Proposition 3.6 is affine, 
so that we can take $C_2=0$. 
Under some integrability condition on $n$, 
quantitative information of $C_2$ 
can be obtained from Lemma 2.6 in \cite{Stannat05} 
combined with Eq. (6) there. 
It will turn out that (\ref{3.18}) is one of basic tools 
in the next section, where more specific cases are discussed. 

A naive guess based on the integration by parts formula (\ref{3.5}) 
would be that ${\rm Sect}(\cE)=\infty$ whenever $a=0$. 
We have not proved this, nor given 
any sufficient condition for the CBCI-process 
not to satisfy the strong sector condition. 
We just present a simple (but very special) 
example of such a CBCI-process. 

\medskip 

\noindent 
{\it Example 3.1}~
Let $b,c>0$ and $n(dy)=c\epsilon_1(dy)$, where 
$\epsilon_1$ is the delta distribution at $1$. 
Then, for each $\delta>0$, 
the CBCI-process with quadruplet $(0,b,n,\delta)$ is ergodic 
and does not satisfy the strong sector condition. 
Indeed, letting $f(x)=\sin 2\pi x$ and $g(x)=\cos 2\pi x$, 
one can observe from (\ref{3.2}) and (\ref{3.5}) that 
$\cE^{\delta}(f,f)=0=\cE^{\delta}(g,g)$ and that 
\begin{eqnarray*}
\cE^{\delta}(f,g) 
& = & 
2\pi c\lg x\int_0^1\cos2\pi(x+y)
(\cos2\pi(x+y)-\cos2\pi x)dy \rg_{\delta}   \\ 
& = & 
\pi c\lg x\rg_{\delta}=\pi c\delta\Phi'(0)=\pi c\delta/b, 
\end{eqnarray*} 
respectively.  

\medskip 

For later use, we close this section 
by giving a lower bound of ${\rm Sect}(\cE)$ 
in a general setting as a refinement of 
the calculation (\ref{3.1}). 
\begin{pr}
Let $\cE$ be as in Proposition 3.1 and 
$f,g\in D(L)$ be such that $\cE(f,f)\cE(g,g)>0$ 
and $\check{\cE}(f,g):=(\cE(f,g)-\cE(g,f))/2>0$. 
Then 
\be 
{\rm Sect}(\cE)^2-1 
\ge 
\left\{
\begin{array}{lll}
\ds{\frac{\check{\cE}(f,g)^2}{\Delta(f,g)}} 
& (\cE(f,f)\cE(g,g)\ne \cE(f,g)\wt{\cE}(f,g), \Delta(f,g)>0) \\ 
\infty & 
(\cE(f,f)\cE(g,g)\ne \cE(f,g)\wt{\cE}(f,g), \Delta(f,g)=0) \\ 
\ds{\frac{\check{\cE}(f,g)}{\wt{\cE}(f,g)}}  & 
(\cE(f,f)\cE(g,g)=\cE(f,g)\wt{\cE}(f,g)),    
\end{array}
\right. 
                                             \label{3.19} 
\ee 
where $\Delta(f,g)=\cE(f,f)\cE(g,g)-\wt{\cE}(f,g)^2\ge 0.$ 
\end{pr}
{\it Proof.}~ 
The proof will be based on 
\be 
{\rm Sect}(\cE)^2\ge 
\sup_{t\in\R}\frac{\cE(f,f+t g)^2}{\cE(f,f)\cE(f+t g,f+t g)}.  
                                             \label{3.20} 
\ee 
Put $U(t)=\cE(f,f+t g)^2/\cE(f+t g,f+t g)$. 
By direct calculations it can be shown 
that $\frac{d}{dt}\log U(t)$ vanishes only for 
\[ 
t=t_0:=\frac{\cE(f,f)(\cE(f,g)-\wt{\cE}(f,g))}
{\cE(f,f)\cE(g,g)-\cE(f,g)\wt{\cE}(f,g)} 
\] 
if $\cE(f,f)\cE(g,g)-\cE(f,g)\wt{\cE}(f,g)\ne 0$, 
whereas $U'(t)$ never vanishes unless $U(t)=0$ 
if $\cE(f,f)\cE(g,g)-\cE(f,g)\wt{\cE}(f,g)= 0$. 
In the former case, the supremum in (\ref{3.20}) is achieved at $t=t_0$ 
and (\ref{3.19}) is obtained for the first two cases 
by calculating $\lim_{t\to t_0}U(t)$. 
In the latter case, 
noting that $\cE(f,g)=\wt{\cE}(f,g)+\check{\cE}(f,g)$, 
we have 
\begin{eqnarray*}
\sup_{t\in\R} \frac{U(t)}{\cE(f,f)} 
& = & 
\lim_{|t|\to\infty} \frac{U(t)}{\cE(f,f)}
\ = \ \frac{\cE(f,g)^2}{\cE(f,f)\cE(g,g)}   \\ 
& = & \frac{\cE(f,g)^2}{\cE(f,g)\wt{\cE}(f,g)} 
\ = \ 1+\frac{\check{\cE}(f,g)}{\wt{\cE}(f,g)},  
\end{eqnarray*} 
which proves (\ref{3.19}) for the third case. 
Lastly, the inequality $\Delta(f,g)\ge 0$ holds in general 
since $\wt{\cE}$ is symmetric and nonnegative definite. 
\qed

\section{Generalized gamma convolutions as stationary distributions}
\setcounter{equation}{0} 

In this section, we apply the results in the previous section 
to the sector constant estimate for a class of CBCI-processes 
whose stationary distributions are GGC's. 
(The general reference for GGC's is \cite{Bondesson}. 
The interested reader is referred also to \cite{SvH04} or \cite{JRY08}.)
The situation, however, is a converse of that in Section 3 
in the following sense. 
We shall be given a priori $q\ge 0$ and 
the L\'evy measure $\Lambda$ of some GGC, 
and then intend to choose $a,b$ and $n$ so that (\ref{3.3}) holds. 
To be more specific, recall that a GGC is 
an infinitely divisible distribution on $\R_+$ with 
L\'evy measure of the form 
\be  
\Lambda_m(dy):=\left(\int e^{-uy}m(du)\right)\frac{dy}{y}, 
                                             \label{4.1} 
\ee 
where $m$ (referred to as the Thorin measure) 
is a measure on $(0,\infty)$ with 
\[ 
\int_{(0,1/2]}|\log u|m(du)
+\int_{(1/2,\infty)}u^{-1}m(du)<\infty. 
\] 
This condition is necessary and sufficient for 
the Laplace exponent 
\be  
\Phi_{q,m}(\lambda)
:=q\lambda+\int \left(1-e^{-\lambda y}\right)\Lambda_m(dy)  
=q\lambda+\int \log\left(1+\frac{\lambda}{u}\right)m(du) 
                                             \label{4.2} 
\ee 
to be finite for all $\lambda>0$. 
(Notice that found in the literature is the condition that 
$\int_{(0,1]}|\log u|m(du)+\int_{(1,\infty)}u^{-1}m(du)<\infty$, 
which allows $m(du)=\one_{(0,1)}(u)du/|\log u|$ inappropriately. 
Here and in what follows, the notation $\one_{E}$ denotes  
the indicator function of a set $E$.) 
We call the distribution having the Laplace exponent (\ref{4.2}) 
the GGC with pair $(q,m)$. 
It must be remembered that 
every gamma distribution has a degenerate Thorin measure. 

In order to study the above mentioned problem, 
we do a heuristic calculation; 
differentiate (\ref{3.3}) with L\'evy measure (\ref{4.1}) to get 
\be  
q+\int \frac{1}{\lambda+u}m(du)
=\frac{1}{\ds{a\lambda+b+\int(1-e^{-\lambda u})\wt{n}(du)}}, 
\qquad \lambda>0. 
                                             \label{4.3} 
\ee 
This equation motivates us to exploit the theory of 
Bernstein functions \cite{SSV}. 
Defining $\cM$ to be the totality of measures $m$ on $(0,\infty)$ 
such that $\int(1+u)^{-1}m(du)<\infty$, 
we recall that every complete Bernstein function $g$ is 
represented uniquely in the form 
\be 
g(\lambda) 
= q\lambda +r+ \int \frac{\lambda}{\lambda+u}m(du), 
\qquad \lambda>0 
                                             \label{4.4} 
\ee 
where $q,r\ge 0$ and $m\in\cM$ (cf. Remark 6.4 in \cite{SSV}), 
and Proposition 7.1 in \cite{SSV} asserts that 
a function $g:(0,\infty)\to\R$ is 
a non-zero complete Bernstein function 
if and only if $g^{\star}(\lambda):=\lambda/g(\lambda)$ 
is a complete Bernstein function. 
With the help of these two facts one can show the next lemma, 
a key to our construction of the desired CBCI-process.  
In what follows we adopt the convention that $1/\infty=0$ 
and set $\overline{m}_{\alpha}=\int u^{\alpha}m(du)$ for $\alpha\in\R$. 
\begin{lm}
Let $q,r\ge 0$ and suppose that 
$m\in\cM$ is non-zero. 
Then there exist uniquely $a,b\ge 0$ and $M\in\cM$ such that 
\be  
q+\frac{r}{\lambda}+\int \frac{1}{\lambda+u}m(du)
=\frac{1}{\ds{a\lambda+b+\int \frac{\lambda}{\lambda+u}M(du)}}, 
\qquad \lambda>0. 
                                             \label{4.5} 
\ee 
Moreover,  
\be 
a=
\left\{
\begin{array}{ll}
0 & (q>0) \\ 
1/(r+\overline{m}_{0}) & (q=0), 
\end{array}
\right. 
\qquad  
\qquad 
b=
\left\{
\begin{array}{ll}
0 & (r>0) \\ 
1/(q+\overline{m}_{-1}) & (r=0)   
\end{array}
\right. 
                                             \label{4.6} 
\ee 
and 
\be 
\overline{M}_{0}=
\left\{
\begin{array}{ll}
1/q-b & (q>0) \\ 
\overline{m}_{1}/(r+\overline{m}_{0})^2-b & 
(q=0, \overline{m}_{0}<\infty) \\ 
\infty & (q=0, \overline{m}_{0}=\infty).   
\end{array}
\right. 
                                             \label{4.7} 
\ee 
\end{lm}
{\it Proof.}~
The first half is immediate from the general facts 
on Bernstein functions previously mentioned. Indeed, 
defining the function $g$ on $(0,\infty)$ by (\ref{4.4}), 
we can find uniquely 
$a,b\ge 0$ and $M\in\cM$ such that 
\[ 
\frac{\lambda}{g(\lambda)}=g^{\star}(\lambda)= 
a\lambda +b+ \int \frac{\lambda}{\lambda+u}M(du), 
\qquad \lambda>0. 
\] 
This is nothing but (\ref{4.5}). 
Most calculations needed to show (\ref{4.6}) and (\ref{4.7}) 
are simple. While letting $\lambda\to\infty$ 
in (\ref{4.5}) yields $a=0$ whenever $q>0$, 
letting $\lambda\downarrow 0$ in (\ref{4.5}) 
gives the value of $b$ in (\ref{4.6}).  In the case $q=0$, 
letting $\lambda\to \infty$ in (\ref{4.5}) multiplied by $\lambda$ 
shows that $a=1/(r+\overline{m}_{0})$. 
With (\ref{4.6}) in mind (\ref{4.7}) 
can be proved in a similar manner because 
$\overline{M}_{0}=\lim_{\lambda\to\infty}
\int \frac{\lambda}{\lambda+u}M(du)$. 
For instance, in the case where $q=0$ and $\overline{m}_{0}<\infty$, 
(\ref{4.5}) and (\ref{4.6}) together yield 
\begin{eqnarray*} 
\int \frac{\lambda}{\lambda+u}M(du)+b 
& = & 
\frac{1}{\ds{\frac{r}{\lambda}+\int \frac{1}{\lambda+u}m(du)}}
-\frac{\lambda}{r+\overline{m}_{0}} \\ 
& = & 
\frac{\ds{\int \frac{\lambda u}{\lambda+u}m(du)}}
{\ds{\left(r+\int \frac{\lambda}{\lambda+u}m(du)\right)
(r+\overline{m}_{0})}}.  
\end{eqnarray*} 
Letting $\lambda\to\infty$ proves (\ref{4.7}) in this case. 
The proof for the other cases are left to the reader. 
\qed

\medskip

\noindent 
Denote by $\cS(m)$ the support of a measure $m$. 
(\ref{4.7}) together with Schwarz's inequality implies that 
$M\equiv 0$ if and only if $q=0=r$ and $m$ is degenerate. 
In such a case, it is understood that $\cS(M)=\emptyset$. 
The next lemma gives bounds of $\inf \cS(M)$ and 
$\sup \cS(M)$ in terms of $q,r$ and $m$. 
\begin{lm}
Let $q,r\ge 0$ and suppose that $m\in\cM$ is non-zero. 
Let $a,b\ge 0$ and $M\in\cM$ be as in Lemma 4.1. 
Define $s_{-}=s_{-}(q,r,m)$ and $s_{+}=s_{+}(q,r,m)$ by 
\[ 
s_{-}=\sup\left\{s<\inf\cS(m):
~s\left(q+\int\frac{m(du)}{u-s}\right)<r \right\}
\] 
and 
\[ 
s_{+}=\inf\left\{s>\sup\cS(m):
~s\left(q-\int\frac{m(du)}{s-u}\right)>r \right\}, 
\] 
respectively. Then the following assertions hold. \\ 
(i) If $a=0=b$, then 
\[ 
\frac{r}{r+\overline{m}_{0}+q\inf\cS(m)}\inf\cS(m) 
\le s_{-}\le \inf\cS(M) 
\] 
and 
\[ 
\sup\cS(M) \le s_{+} 
\left\{
\begin{array}{ll}
\le \sup\cS(m)+(r+\overline{m}_{0})q^{-1} & (q>0) \\ 
= \infty & (q=0).  
\end{array} 
\right. 
\] 
(ii) If $a>0$ and $b=0$, then 
\[ 
\frac{r}{r+\overline{m}_{0}} \inf\cS(m) 
\le s_{-}\le \inf\cS(M)\le \sup\cS(M) \le \sup\cS(m).  
\] 
(iii) If $a=0$ and $b>0$, then 
\[ 
\inf\cS(m)\le \inf\cS(M)\le \sup\cS(M) \le s_{+}
\left\{
\begin{array}{ll}
\le \sup\cS(m)+\overline{m}_{0}q^{-1} & (q>0) \\ 
= \infty & (q=0).
\end{array} 
\right. 
\] 
(iv) If $a,b>0$, then 
$\cS(M)$ is contained in the closed interval from $\inf\cS(m)$ 
to $\sup\cS(m)$. 
\end{lm}
Before going to the proof, it is worth noting that 
the analytic extension of (\ref{4.5}) can be regarded as 
a relation between Stieltjes transforms of $m$ and $M$. 
Introducing the notation 
$G_m(z)=\int(z-u)^{-1}m(du)$ for $m\in \cM$, 
we deduce from (\ref{4.5}) 
\be 
G_m(z)-q+\frac{r}{z}
=\frac{1}{az-b-zG_M(z)}, 
\qquad z\in\C\setminus [\inf\cS(m),\sup\cS(m)]. 
                                         \label{4.8}
\ee 
{\it Proof of Lemma 4.2.}~
We shall employ the following fact. 
(See e.g. Theorem A.6 in \cite{KN}.) 
For any $s_1<s_2$, a measure 
$M$ is supported on $[s_1,s_2]$ if and only if 
$G_M$ is holomorphic on $\C\setminus [s_1,s_2]$, 
negative on the interval $(-\infty,s_1)$ 
and positive on the interval $(s_2,\infty)$. 

To show (i) we assume, in addition to $a=b=0$, that
$0<\inf \cS(m)\le\sup\cS(m)<\infty$. 
By (\ref{4.8}) we have 
\[ 
G_M(z)
= \frac{1}{\ds{z\left(q-\int\frac{m(du)}{z-u}\right)-r}}
= \frac{-1}{\ds{r-z\left(q+\int\frac{m(du)}{u-z}\right)}},  
\] 
from which it is easily verified that $G_M$ 
is holomorphic on $\C\setminus [s_{-},s_{+}]$, 
positive for all $z>s_{+}$ and negative for all $z<s_{-}$. 
So, $s_{-}\le \inf\cS(M)\le \sup\cS(M) \le s_{+}$ 
for the abovementioned reason.  
To prove the estimates for $s_{-}$ and $s_{+}$, 
we may assume further that $\overline{m}_{0}<\infty$. 
We then have $q>0$ because of (\ref{4.6}). 
It is clear that a lower estimate for $s_{-}$ is given by 
the smallest solution, say $s_{-}'=s_{-}'(q,r,m)$, 
to the quadratic equation in $s$ 
\[ 
s\left(q+\frac{\overline{m}_{0}}{\inf\cS(m)-s}\right)=r, 
\] 
namely 
\begin{eqnarray} 
s_{-}
& \ge &  s_{-}'                                        \nonumber \\ 
& = & 
\frac{(r+\overline{m}_{0}+q\inf\cS(m))
-\sqrt{(r+\overline{m}_{0}+q\inf\cS(m))^2-4qr \inf\cS(m)}}{2q} 
                                                    \label{4.9} \\ 
& \ge & 
\frac{r}{r+\overline{m}_{0}+q\inf\cS(m)}\inf\cS(m).    \nonumber 
\end{eqnarray} 
An upper estimate for $s_{+}$ is given by 
the largest solution, say $s_{+}'=s_{+}'(q,r,m)$, to 
\[ 
s\left( q+\frac{\overline{m}_{0}}{\sup\cS(m)-s}\right)=r. 
\] 
Thus $s_{+}$ is dominated by 
\begin{eqnarray*} 
s_{+}'
& = & 
\frac{(r+\overline{m}_{0}+q\sup\cS(m))
+\sqrt{(r+\overline{m}_{0}+q\sup\cS(m))^2-4qr \sup\cS(m)}}{2q} \\
& \le & 
\sup\cS(m)+\frac{r+\overline{m}_{0}}{q}. 
\end{eqnarray*} 
Also, it is obvious that $s_{+}(=\inf\emptyset)=\infty$ for $q=0$. 

Next we prove (ii), assuming that $a>0$ and $b=0$. 
By the former $q=0$ and $\overline{m}_{0}<\infty$. 
Therefore, again by (\ref{4.8}) 
\[ 
G_M(z)
= \frac{\ds{\int\frac{u}{z-u}m(du)}}
{\ds{\left(r-z\int\frac{m(du)}{u-z}\right)(r+\overline{m}_{0})}}. 
\] 
This allows proceeding along the same lines as the proof of (i). 
The proof of (iii) follows very closely the proof of (ii). 
So, these proofs are omitted. 
It remains to prove (iv). 
To this end, assume that $a,b>0$ and observe 
from (\ref{4.6}) and (\ref{4.5}) that 
\[ 
G_M(z)= 
\frac{\ds{\overline{m}_{-1}\int\frac{1}{z-u}m(du)
-\overline{m}_{0}\int\frac{1}{u(z-u)}m(du)}}
{\ds{\overline{m}_{0}\overline{m}_{-1}\int\frac{1}{z-u}m(du)}}, 
\] 
in which both $\overline{m}_{0}$ and $\overline{m}_{-1}$ 
are positive and finite. 
This shows the analyticity of $G_M$ 
on $\C\setminus [\inf\cS(m),\sup\cS(m)]$. 
We may assume additionally that $m$ is nondegenerate, 
for otherwise the assertion is trivial since $M\equiv 0$. 
Then the positivity of the numerator 
on the right side for any 
$z\in\R\setminus [\inf\cS(m),\sup\cS(m)]$ 
follows from 
\begin{eqnarray*}
\lefteqn{\overline{m}_{-1}\int\frac{1}{z-u}m(du)
-\overline{m}_{0}\int\frac{1}{u(z-u)}m(du)} \\ 
& = & 
\int\int\left(\frac{1}{v(z-u)}-\frac{1}{v(z-v)}\right)m(du)m(dv) \\ 
& = & 
\int\frac{u}{z-u}m(du)\int\frac{1}{v(z-v)}m(dv)
-\left(\int\frac{1}{z-u}m(du)\right)^2>0. 
\end{eqnarray*} 
These properties together prove (iv). 
\qed 

\medskip

We now present the main result of this section, 
which concerns the construction of the CBCI-process 
having a given GGC as a stationary distribution. 
In other words, the problem (II) addressed in Introduction is solved. 
Simultaneously, the sector constant estimate will be obtained 
as a solution to the problem (III). 
Note that every Thorin measure belongs to $\cM$. 
\begin{th}
Let $q\ge 0$ and suppose that $m$ is a non-zero Thorin measure. 
Let $a,b\ge 0$ and $M\in\cM$ be as in Lemma 4.1 with $r=0$. 
Then the GGC with pair $(q,m)$ is a unique 
stationary distribution of the CBCI-process 
with quadruplet $(a,b,n,1)$, where $n$ is a measure on $(0,\infty)$ 
defined to be 
\be  
n(dy) =dy \int u^2e^{-uy}M(du). 
                                             \label{4.10} 
\ee 
If, in addition, $q=0, \overline{m}_{1}<\infty$ 
and $\inf\cS(m)>0$, then $\overline{m}_{0}<\infty$, $\inf\cS(M)>0$ and 
the CBCI-process with quadruplet $(a,b,n,\delta)$ satisfies 
\be  
{\rm Sect}(\cE^{\delta})-1 
\le 
\sqrt{\left(\frac{\overline{m}_{1}}
{\overline{m}_{0}}-\inf\cS(m)\right)\frac{2}{\inf\cS(M)}} 
\le 
\sqrt{2\left(\frac{\overline{m}_{1}}
{\overline{m}_{0}\cdot\inf\cS(m)}-1\right)} 
                                             \label{4.11} 
\ee 
for any $\delta>0$. 
\end{th}
{\it Proof}~ 
Firstly, we claim $\int\min\{y,y^2\}n(dy)<\infty$. 
To see this, note that by (\ref{4.10}) and Fubini's theorem 
\begin{eqnarray}
\int\min\{y,y^2\}n(dy) 
& = & \int M(du)u^2\int \min\{y,y^2\}e^{-uy}dy \nonumber \\ 
& = & \int M(du)\int y\min\{1,y/u\}e^{-y}dy \nonumber \\ 
& = & \int M(du)
\left(\int_0^u (y^2/u)e^{-y}dy+\int_u^{\infty} ye^{-y}dy \right).  
                                             \label{4.12} 
\end{eqnarray} 
It is elementary to verify that the integrand 
in the last expression is bounded above (and below)   
by a positive constant times $1/(1+u)$. 
Therefore it follows from $M\in\cM$ that 
$\int\min\{y,y^2\}n(dy)$ is finite. 
Applying Fubini's theorem, we see that 
$\int(1-e^{-\lambda y})\wt{n}(dy) 
=\lambda\int(\lambda+u)^{-1}M(du)$ for $\lambda>0$. 
Plugging this into (\ref{4.5}) with $r=0$ yields (\ref{4.3}). 
By integrating it 
\[ 
q\lambda + \int \log\left(1+\frac{\lambda}{u}\right)m(du) 
=\int_{0}^{\lambda} 
\frac{du}{\ds{au+b+\int(1-e^{-uy})\wt{n}(dy)}}, 
\qquad \lambda\ge 0. 
\] 
In view of (\ref{3.3}) and (\ref{4.2}), 
this proves the first half. 

For the proof of the last half, we assume additionally that 
$q=0$, $\overline{m}_{1}<\infty$ and $\inf\cS(m)>0$. 
The last two conditions together imply that 
$\overline{m}_{0}, \overline{m}_{-1}<\infty$ 
since $m$ is a Thorin measure.  
Hence $a,b>0$ by (\ref{4.6}). 
Also, by (\ref{4.10}) and (\ref{4.7}) 
$\int yn(dy)=\overline{M}_{0} 
=\overline{m}_{1}/\overline{m}_{0}^2-1/\overline{m}_{-1}$, 
which is finite. 
We now apply Proposition 3.6 
to the L\'evy measure (\ref{4.1}) or equivalently 
to the function $\varphi(y):=\int e^{-uy} m(du)$. 
It is easily observed that $V(y)=-\log \varphi(y)$ satisfies 
\[  
V'(0)-V'(y) 
=  
\frac{\ds{\int um(du)}}{\ds{\int m(du)}}
-\frac{\ds{\int e^{-uy}um(du)}}{\ds{\int e^{-uy}m(du)}}  
\le  
\frac{\overline{m}_{1}}{\overline{m}_{0}}-\inf\cS(m)    
\] 
for all $y>0$. 
On the other hand, it follows from (\ref{4.10}) that 
\[ 
\left\Vert \frac{d\wt{n}}{dn}\right\Vert_{\infty}
\le 
\frac{1}{\inf\cS(M)} 
\le 
\frac{1}{\inf\cS(m)}, 
\] 
where the last inequality is implied by Lemma 4.2 (iv).  
Therefore, (\ref{4.11}) is deduced from (\ref{3.18}). 
\qed 

\medskip 

\noindent 
Notice that the upper estimates (\ref{4.11}) are effective in the sense that 
the most right side vanishes only in the reversible case, 
namely when $m$ is degenerate. 
The symmetry found in the statement of Lemma 4.1 
allows one to show 
a converse of Theorem 4.3. 
Roughly speaking, we will see below that 
every ergodic CBCI-process with $n$ having 
a (non-zero) completely monotone density 
has some GGC as a (non-reversible) stationary distribution.  
\begin{th}
Let $a,b\ge 0$ and suppose that a non-zero $M\in\cM$ is given. 
Define $n$ and $\Phi$ by (\ref{4.10}) and 
the first equality in (\ref{3.3}), respectively. 
Then $\int\min\{y,y^2\}n(dy)$ is finite and 
the following assertions hold true. \\ 
(i) $\Phi(1)<\infty$ if and only if $b>0$ or 
\[ 
\int_{0}^1\left(u\ds{\int \frac{M(dv)}{u+v}}\right)^{-1}du<\infty. 
\] 
(ii) If $\Phi(1)<\infty$, then 
the unique stationary distribution of 
the CBCI-process with quadruplet $(a,b,n,1)$ 
is a GGC with some pair $(q,m)$ such that 
\be 
q=
\left\{
\begin{array}{ll}
0 & (a>0) \\ 
1/(b+\overline{M}_{0}) & (a=0), 
\end{array}
\right. 
\quad 
\overline{m}_{0}=
\left\{
\begin{array}{ll}
1/a & (a>0) \\ 
\overline{M}_{1}/(b+\overline{M}_{0})^2 & 
(a=0, \overline{M}_{0}<\infty) \\ 
\infty & (a=0, \overline{M}_{0}=\infty).   
\end{array}
\right. 
                                             \label{4.13} 
\ee 
(iii) If $a,b>0, \overline{M}_{0}<\infty$ 
and $\inf\cS(M)>0$, then for any $\delta>0$ 
the CBCI-process with quadruplet $(a,b,n,\delta)$ is ergodic and 
$\ds{\left({\rm Sect}(\cE^{\delta})-1\right)^2}$ is dominated by 
\be 
\frac{\sqrt{\left(a\inf\cS(M)-b-\overline{M}_{0}\right)^2
+4a\overline{M}_{0}\inf\cS(M)}
-\left(a\inf\cS(M)-b-\overline{M}_{0}\right)}{a\inf\cS(M)}. 
                                             \label{4.14} 
\ee 
\end{th}
{\it Proof.}~
By virtue of (\ref{4.12}), 
that $\int\min\{y,y^2\}n(dy)<\infty$ is ensured by $M\in\cM$. 
Putting $g_M(u)=\int(u+v)^{-1}M(dv)$, note that 
$\Phi(\lambda)= 
\int_0^{\lambda}(au+b+ug_M(u))^{-1}du.$ 
For the proof of (i), it is sufficient to show that $b=0$ and 
$\int_0^{1}(ug_M(u))^{-1}du =\infty$ together imply $\Phi(1)=\infty$. 
If $\overline{M}_{-1}<\infty$, 
this holds true since 
$\Phi(1)\ge\int_0^{1}(au+u\overline{M}_{-1})^{-1}du=\infty$. 
If $\overline{M}_{-1}=\infty$, 
there exists $v\in(0,1)$ such that 
$a<g_M(u)$ for any $u\in(0,v]$ 
and hence $\Phi(1)\ge \int_0^{v}(2ug_M(u))^{-1}du=\infty$. 
These observations prove (i). 

To show (ii) we interchange in Lemma 4.1 
the roles of $(a,b,M)$ and $(q,r,m)$ to get 
\be  
a+\frac{b}{\lambda}+\int \frac{1}{\lambda+u}M(du)
=\frac{1}{\ds{q\lambda+r+\int \frac{\lambda}{\lambda+u}m(du)}}, 
\qquad \lambda>0 
                                             \label{4.15} 
\ee 
for some $q,r\ge 0$ and $m\in\cM$. 
Here, according to (\ref{4.6}), $q$ is given by (\ref{4.13}), 
$r=0$ if $b>0$, and $r=1/(a+ \overline{M}_{-1})$ if $b=0$. 
But it follows from $\Phi(1)<\infty$ and (i) that 
$\overline{M}_{-1}=\infty$ whenever $b=0$. 
As a result $r=0$, and therefore the formula (\ref{4.13}) 
for $\overline{m}_{0}$ precisely corresponds to (\ref{4.7}). 
By (\ref{4.15}) and Fubini's theorem 
one can show that $\Phi(\lambda)=\Phi_{q,m}(\lambda)$. 
Since $\Phi(1)<\infty$, this proves not only that 
$m$ is a Thorin measure but also the assertion (ii). 

It remains to show the sector constant estimate 
under the stronger assumptions that $a,b>0, \overline{M}_{0}<\infty$ 
and $\inf\cS(M)>0$. This should be derived from (\ref{4.11}) 
once $\overline{m}_{1}/\overline{m}_{0}$ and 
$\inf \cS(m)$ are estimated in terms of $a,b$ and $M$. 
For this purpose, combine $\overline{m}_{0}=1/a$ with 
$\overline{M}_{0}=\overline{m}_{1}/{\overline{m}_{0}}^2-b$ 
to get $\overline{m}_{1}/\overline{m}_{0}=(\overline{M}_{0}+b)/a$. 
Also, Lemma 4.2 (i) and (\ref{4.9}) 
with $(q,r,M)$ being interchanged with $(a,b,m)$ are 
applied to derive 
\begin{eqnarray*} 
\inf\cS(M) 
& \ge & s_{-}(a,b,M) \ \ge \ s_{-}'(a,b,M) \\ 
& = & 
\frac{\left(a\inf\cS(M)+b+\overline{M}_{0}\right)
-\sqrt{\left(a\inf\cS(M)+b+\overline{M}_{0}\right)^2-4ab\inf\cS(M)}}{2a}. 
\end{eqnarray*} 
These calculations together yield 
\begin{eqnarray*} 
\lefteqn{\frac{\overline{m}_{1}}{\overline{m}_{0}}-\inf\cS(m)}     \\
& \le & 
\frac{\left(b+\overline{M}_{0}\right)-a\inf\cS(M)
+\sqrt{\left(a\inf\cS(M)+b+\overline{M}_{0}\right)^2-4ab\inf\cS(M)}}{2a} \\ 
& = & 
\frac{\sqrt{\left(a\inf\cS(M)-b-\overline{M}_{0}\right)^2
+4a\overline{M}_{0}\inf\cS(M)}
-\left(a\inf\cS(M)-b-\overline{M}_{0}\right)}{2a}.  
\end{eqnarray*} 
Therefore,  the bound (\ref{4.14}) is deduced from 
the first inequality in (\ref{4.11}). 
\qed 

\medskip 

To check, consider the case discussed in Example 2.1 (ii). 
It is easily seen that the Laplace exponent (\ref{2.24}) 
is that of the GGC with $q=0$ and Thorin measure 
\be 
m(du)=\frac{1}{\Gamma(\alpha)\Gamma(1-\alpha)}\cdot
\frac{\one_{(\kappa,\infty)}(u)}{(u-\kappa)^{\alpha}}du. 
                                             \label{4.16} 
\ee 
Accordingly, $\overline{m}_{0}=\infty$ and 
$\overline{m}_{-1}=\kappa^{-\alpha}$. 
So, (\ref{4.6}) gives $a=0$ and $b=\kappa^{\alpha}$ consistently. 
Although Lemma 4.1 itself does not give 
any explicit form of $M$, 
we can identify it with 
\be 
M(du)=\frac{u^{-1}(u-\kappa)^{\alpha}}{\Gamma(\alpha)\Gamma(1-\alpha)}
\cdot\one_{(\kappa,\infty)}(u)du, 
                                             \label{4.17} 
\ee 
thanks to (\ref{2.23}). Therefore, in this case 
$\cS(M)=\cS(m)=(\kappa,\infty)$. 
In the next section we will be provided with 
some procedure to derive a formulae for $M$ including (\ref{4.17}) 
via (\ref{4.5}) and with further examples as well. 

As for the lower bound of the sector constant of 
CBCI-processes discussed in Theorem 4.3, 
one can show the next result as an application of Proposition 3.7. 
\begin{th}
Suppose that $m$ is a non-zero measure on $(0,\infty)$ 
with $\inf\cS(m)>0$ and $\overline{m}_{0}<\infty$. 
Let $a,b>0$ and $M\in\cM$ be as in Lemma 4.1 with $q=0=r$ 
and $n$ be given by (\ref{4.10}). 
Then the bilinear form $\cE$ associated with 
the CBCI-process with quadruplet $(a,b,n,1)$ satisfies 
\[  
{\rm Sect}(\cE)^2-1 
\ge 
\frac{\left(\overline{m}_{-1}\overline{m}_{-3}-
\overline{m}_{-2}^2\right)^2}
{2\overline{m}_{-1}\overline{m}_{-2}^2\overline{m}_{-3}
+4\overline{m}_{-1}^2\overline{m}_{-2}^3
+12\overline{m}_{-1}^2\overline{m}_{-2}\overline{m}_{-4}
-9\overline{m}_{-1}^2\overline{m}_{-3}^2-\overline{m}_{-2}^4}. 
\] 
\end{th}
{\it Proof}~ 
By the assumption $m$ is a Thorin measure and 
the GGC with pair $(0,m)$ has an exponential integrability. 
Indeed, its Laplace exponent 
$\Phi_{0,m}(\lambda)$ can extend 
real analytically to $\lambda>-\inf \cS(m)$. In particular, 
we have the finite moments given by 
\[ 
\lg x^k\rg = (-1)^k
\left.\frac{d^k}{d\lambda^k}e^{-\Phi_{0,m}(\lambda)}\right|_{\lambda=0}
\] 
for $k=1,2,\ldots.$ Since 
$\Phi_{0,m}^{(k)}(0)=(-1)^{k-1}(k-1)!\overline{m}_{-k}$, 
it follows that 
\be 
\lg x\rg = \overline{m}_{-1}, 
\ 
\lg x^2\rg = \overline{m}_{-1}^2+\overline{m}_{-2} 
\ \mbox{and} \ 
\lg x^3\rg = 2\overline{m}_{-3}
+3\overline{m}_{-1}\overline{m}_{-2}+\overline{m}_{-1}^3. 
                                             \label{4.18} 
\ee 
We are going to apply (\ref{3.19}) by 
taking $f(x)=x$ and $g(x)=x^2$, for which 
$\Delta(f,g)>0$ holds by virtue of 
(\ref{3.2}) and Schwarz's inequality together. 
Explicitly, our task is to show for some common $C>0$ that 
\be 
\check{\cE}(f,g)^2
=C(\overline{m}_{-1}\overline{m}_{-3}-\overline{m}_{-2}^2)^2 
                                             \label{4.19} 
\ee 
and that $\Delta(f,g)={\cE}(f,f){\cE}(g,g)-\wt{\cE}(f,g)^2$ equals 
\be 
C(2\overline{m}_{-1}\overline{m}_{-2}^2\overline{m}_{-3}
+4\overline{m}_{-1}^2\overline{m}_{-2}^3
+12\overline{m}_{-1}^2\overline{m}_{-2}\overline{m}_{-4}
-9\overline{m}_{-1}^2\overline{m}_{-3}^2-\overline{m}_{-2}^4). 
                                             \label{4.20} 
\ee 
Noting that $\int n(dy)y^k=k!\overline{M}_{1-k}$ $(k=1,2,\ldots)$ by 
(\ref{4.10}), observe from (an extension of) (\ref{3.5}) and (\ref{3.2}) 
that 
\be 
\check{\cE}(f,g)
=\frac{1}{2}\lg x\int \wt{n}(dy)(-y^2)\rg
=-\frac{1}{6}\lg x\rg \int n(dy)y^3=- \lg x\rg \overline{M}_{-2}, 
                                             \label{4.21} 
\ee 
\[ 
{\cE}(f,f)
=\lg x\rg\left(a+\frac{1}{2}\int n(dy)y^2\right) 
=\lg x\rg\left(a+\overline{M}_{-1}\right), 
\] 
\begin{eqnarray*}
\wt{\cE}(f,g) 
& = & 
2\lg x^2\rg\left(a+\frac{1}{2}\int n(dy)y^2\right)
+\frac{1}{2}\lg x\rg \int n(dy)y^3  \\ 
& = & 
2\lg x^2\rg\left(a+\overline{M}_{-1}\right)
+3\lg x\rg \overline{M}_{-2} 
\end{eqnarray*} 
and 
\begin{eqnarray*} 
\cE(g,g) 
& = & 
4\lg x^3\rg\left(a+\frac{1}{2}\int n(dy)y^2\right) 
+2\lg x^2\rg \int n(dy)y^3
+\frac{1}{2}\lg x\rg \int n(dy)y^4   \\ 
& = & 
4\lg x^3\rg\left(a+\overline{M}_{-1}\right) 
+12\lg x^2\rg\overline{M}_{-2}
+12\lg x\rg \overline{M}_{-3}.   
\end{eqnarray*} 
The last three equalities together yield 
\begin{eqnarray} 
\Delta(f,g) 
& = & 
4\left(\lg x\rg\lg x^3\rg-\lg x^2\rg^2\right)
\left(a+\overline{M}_{-1}\right)^2 
                                          \nonumber \\ 
&  & 
+12\lg x\rg^2 \left(a+\overline{M}_{-1}\right)\overline{M}_{-3} 
-9\lg x\rg^2 {\overline{M}_{-2}}^2.   
                                             \label{4.22} 
\end{eqnarray} 
We shall calculate 
$\overline{M}_{-1}$, $\overline{M}_{-2}$ and 
$\overline{M}_{-3}$. The results are 
\be 
\overline{M}_{-1} = 
\frac{\overline{m}_{0}\overline{m}_{-2}
-\overline{m}_{-1}^2}{\overline{m}_{0}\overline{m}_{-1}^2}, 
\ \  
\overline{M}_{-2} = 
\frac{\overline{m}_{-1}\overline{m}_{-3}
-\overline{m}_{-2}^2}{\overline{m}_{-1}^3} 
                                             \label{4.23} 
\ee 
and 
\be 
\overline{M}_{-3}
= 
\frac{\overline{m}_{-1}^2\overline{m}_{-4}
-2\overline{m}_{-1}\overline{m}_{-2}\overline{m}_{-3}
+\overline{m}_{-2}^3}{\overline{m}_{-1}^4},  
                                             \label{4.24} 
\ee 
whose proof are postponed for a while. 
Then (\ref{4.19}) with $C=1/\overline{m}_{-1}^4$ 
follows immediately by plugging (\ref{4.18}) 
and (\ref{4.23}) into (\ref{4.21}). 
In principle, (\ref{4.20}) with the same $C$ 
can be obtained similarly though the calculation is tedious. 
To carry this out, observe that  
$a+\overline{M}_{-1}=\overline{m}_{-2}/\overline{m}_{-1}^2$ 
by $a=1/\overline{m}_{0}$ and that (\ref{4.18}) give  
\[ 
\lg x\rg\lg x^3\rg-\lg x^2\rg^2
=\overline{m}_{-1}^2\overline{m}_{-2}
+2\overline{m}_{-1}\overline{m}_{-3}-\overline{m}_{-2}^2. 
\] 
(\ref{4.20}) with $C=1/\overline{m}_{-1}^{4}$ 
will now be derived in a fairly straightforward way. 

It remains to prove (\ref{4.23}) and (\ref{4.24}). 
We exploit a variant of the identity used 
to find $\overline{M}_{0}$ in the proof of Lemma 4.1 (with $q=0=r$): 
\be 
\int \frac{1}{\lambda+u}M(du)
= 
\frac{\ds{\int \frac{u}{\lambda+u}m(du)}}
{\ds{\overline{m}_{0}}\int \frac{\lambda}{\lambda+u}m(du)}
-\frac{1}{\overline{m}_{-1} \lambda}, 
\qquad \lambda>0. 
                                             \label{4.25} 
\ee 
Letting $\lambda\downarrow 0$, we get 
(\ref{4.23}) for $\overline{M}_{-1}$ 
with the help of L'Hospital's rule. 
By differentiating (\ref{4.25}) 
\begin{eqnarray} 
-\int \frac{1}{(\lambda+u)^2}M(du) 
& = & 
\frac{1}{\lambda^2}\cdot
\frac{\ds{-\int \frac{um(du)}{(\lambda+u)^2}\int \frac{m(dv)}{v}
+\left(\int \frac{m(du)}{\lambda+u}\right)^2}}
{\ds{\overline{m}_{-1}\left(\int \frac{m(du)}{\lambda+u}\right)^2}} 
                                       \nonumber \\ 
& = & 
-\frac{\ds{\int\int \frac{(u-v)^2}
{(\lambda+u)^2uv(\lambda+v)^2}m(du)m(dv)}}
{2\ds{\overline{m}_{-1}\left(\int \frac{m(du)}{\lambda+u}\right)^2}},   
                                             \label{4.26} 
\end{eqnarray} 
where the symmetry of $m(du)m(dv)$ applies to show 
the last equality. Letting $\lambda\downarrow 0$ 
leads to (\ref{4.23}) for $\overline{M}_{-2}$. 
Finally, differentiating the both sides of (\ref{4.26}) 
multiplied by the square of $\int({\lambda+u})^{-1}m(du)$ 
and then letting $\lambda\downarrow 0$, 
we see without difficulty that 
\[ 
2\overline{M}_{-3}\overline{m}_{-1}^2
+2\overline{M}_{-2}\overline{m}_{-1}\overline{m}_{-2}
= 
\frac{2(\overline{m}_{-1}\overline{m}_{-4}
-\overline{m}_{-2}\overline{m}_{-3})}{\overline{m}_{-1}}. 
\] 
This combined with (\ref{4.23}) proves (\ref{4.24}). 
The proof of Theorem 4.5 is complete. 
\qed

\section{Further discussions and related topics}
\setcounter{equation}{0} 

In this section, most calculations will be 
based on the equations (\ref{4.5}) and (\ref{4.8}) with $r=0$, 
and therefore we take $r=0$ without explicit mention 
(except in the statements of Propositions 5.2 and 5.3 below). 
The first two subsections are devoted to further studies  
of the correspondence between $(a,b,M)$ and $(q,m)$ 
under some special circumstances. 
In the final subsection, the basic relation (\ref{4.5}) 
will be discussed also in connection with 
certain topics in noncommutative probability theory. 

\subsection{Discrete Thorin measures} 

This subsection concerns  
the correspondence between GGC's and ergodic CBCI processes 
under the condition that $m$ (or $M$) is discrete. 
As the simplest example of such GGC's (other than gamma distributions) 
we first discuss the case of convolutions of two gamma distributions 
in rather detail. Given $x\in\R$, 
denote by $\epsilon_x$ the delta distribution at $x$. \\ 
{\it Example 5.1.}~
(i) 
Let $\gamma_1,\gamma_2,\lambda_1,\lambda_2$ be positive constants. 
Consider $m:=\gamma_1\epsilon_{\lambda_1}+\gamma_2\epsilon_{\lambda_2}$ 
as the Thorin measure. 
According to (\ref{4.6}) with $q=0$, 
$a$ and $b$ are chosen as
\be  
a=\frac{1}{\gamma_1+\gamma_2} 
\quad \mbox{and} \quad 
b=\frac{1}{\lambda_1^{-1}\gamma_1+\lambda_2^{-1}\gamma_2}, 
                                             \label{5.1} 
\ee
respectively. Our ansatz here is that 
the measure $M$ satisfying (\ref{4.5}) 
with $q=0$ is of the form 
$M=c\epsilon_{\kappa}$ for some $c,\kappa>0$, 
which are to be determined. 
The measure $n$ in (\ref{4.10}) is then given by 
$n(dy)=c\kappa^2e^{-\kappa y}dy$ and 
the equation (\ref{4.5}) reads 
\be  
\frac{\gamma_1}{\lambda+\lambda_1}
+\frac{\gamma_2}{\lambda+\lambda_2}
=\frac{1}{\ds{a\lambda+b+\frac{c\lambda}{\lambda+\kappa}}}, 
\qquad \lambda\ge 0. 
                                             \label{5.2} 
\ee 
It is not difficult to see that 
the above requirement is fulfilled by setting 
\be 
\kappa:=\frac{a}{b}\lambda_1\lambda_2
=\frac{\lambda_2\gamma_1+\lambda_1\gamma_2}{\gamma_1+\gamma_2} 
                                             \label{5.3} 
\ee  
and 
\be  
c:=a(\lambda_1+\lambda_2-\kappa)-b 
=\frac{\gamma_1\gamma_2(\lambda_1-\lambda_2)^2}
{(\gamma_1+\gamma_2)^2(\lambda_2\gamma_1+\lambda_1\gamma_2)},  
                                             \label{5.4} 
\ee 
which vanishes for degenerate $m$ with $\lambda_1=\lambda_2$  
in accordance with the comments 
in the paragraph preceding to Lemma 4.2. 
Lastly, assuming $\lambda_1\le \lambda_2$ 
and letting $\delta>0$, we have 
by the first bound in (\ref{4.11}) combined with (\ref{5.3})
\be 
{\rm Sect}(\cE^{\delta})-1 
\le 
\sqrt{\left(\frac{\gamma_1\lambda_1+\gamma_2\lambda_2}
{\gamma_1+\gamma_2}-\lambda_1\right)\frac{2}{\kappa}}
=
\sqrt{\frac{2(\lambda_2-\lambda_1)\gamma_2}
{\lambda_2\gamma_1+\lambda_1\gamma_2}} 
                                             \label{5.5} 
\ee 
for the CBCI-process with quadruplet $(a,b,n,\delta)$. 

\noindent 
(ii) One can reverse the above procedure 
to construct explicitly the L\'evy density of 
the stationary distribution of the CBCI-process 
with $n$ being of the form $n(dy)=c\kappa^2e^{-\kappa y}dy$. 
Indeed, given $a,b,c,\kappa>0$, we find 
$\lambda_1,\lambda_2,\gamma_1,\gamma_2$ satisfying (\ref{5.2}) 
in the following manner. In view of the first equalities of 
(\ref{5.3}) and (\ref{5.4}), 
the required $\lambda_i$'s must solve the equation 
$p(\lambda):=a\lambda^2-(a\kappa+b+c)\lambda+b\kappa=0$. 
This leads to 
\[ 
\lambda_1=\frac{(a\kappa+b+c)-\sqrt{D}}{2a}
\quad \mbox{and} \quad 
\lambda_2=\frac{(a\kappa+b+c)+\sqrt{D}}{2a}, 
\] 
where $D=(a\kappa+b+c)^2-4ab\kappa
=(a\kappa-b-c)^2+4ac\kappa>0$. 
Note that $\lambda_1<\kappa<\lambda_2$ since $p(\kappa)<0$. 
Moreover, $\gamma_i$'s are determined by 
\[  
\gamma_1=\frac{\kappa-\lambda_1}{a(\lambda_2-\lambda_1)}
\quad \mbox{and} \quad 
\gamma_2=\frac{\lambda_2-\kappa}{a(\lambda_2-\lambda_1)},  
\]  
for which (\ref{5.1}) are easily checked to hold. 
Consequently, the L\'evy density (\ref{4.1}) equals 
$(\gamma_1 e^{-\lambda_1 y}+\gamma_2e^{-\lambda_2 y})/y$. 
In addition, one can derive from (\ref{5.5}) 
the sector constant estimate described 
in terms of $a,b,c$ and $\kappa$ by noting that 
\[ 
\frac{2(\lambda_2-\lambda_1)\gamma_2}
{\lambda_2\gamma_1+\lambda_1\gamma_2}
=\frac{2(\lambda_2-\kappa)}{\kappa}
=\frac{\sqrt{D}-(a\kappa-b-c)}{a\kappa}, 
\] 
which clearly corresponds to (\ref{4.14}). 

\smallskip 

\noindent 
{\it Example 5.2.}~
Given $q>0$ and a degenerate Thorin measure $m$, 
we have a similar situation to Example 5.1 
except the sector constant estimate, 
which is not available since $a=0$ by (\ref{4.6}). 
Indeed, for $m=\gamma_1\epsilon_{\lambda_1}$ 
with $\gamma_1,\lambda_1>0$ being fixed arbitrarily, 
$b=1/(q+\gamma_1/\lambda_1)$ by (\ref{4.6}) and 
it is easily verified that 
the measure $M$ satisfying (\ref{4.5}) is given by  
\[ 
M=\left(\frac{1}{q}-\frac{1}{q+\gamma_1/\lambda_1}\right)
\epsilon_{\lambda_1+\gamma_1/q}. 
\] 
Conversely, if we are given $b>0$ and $M=c\epsilon_{\kappa}$ 
with $c,\kappa>0$, then the measure $m$ determined by (\ref{4.5}) 
with $q=1/(b+c)$ and $a=0$ is shown to be 
\[ 
m=\frac{c\kappa}{(b+c)^2}\epsilon_{b\kappa/(b+c)}. 
\] 
Alternatively, this can be derived directly from 
(\ref{3.4}) by noting that 
$\wt{n}^{*N}(dy)=(c\kappa)^Ny^{N-1}e^{-\kappa y}dy/(N-1)!$ 
for $N=1,2,\ldots.$

\medskip 

Apart from explicit expressions, 
the above examples are generalized as follows. 
\begin{pr}
Let $l\ge 1$ be a fixed integer. \\ 
(i) For every discrete measure  
\be 
m=\gamma_1\epsilon_{\lambda_1}
+\cdots+\gamma_{l+1}\epsilon_{\lambda_{l+1}} 
\quad \mbox{with} \quad  
\gamma_i>0 \ \mbox{and} \ 0<\lambda_1<\cdots<\lambda_{l+1}, 
                                             \label{5.6} 
\ee 
there exists a unique measure $M$ of the form 
\be 
M=c_1\epsilon_{\kappa_1}+\cdots+c_{l}\epsilon_{\kappa_{l}} 
\quad \mbox{with} \quad  
c_i>0 \ \mbox{and} \ 0<\kappa_1<\cdots<\kappa_{l} 
                                             \label{5.7} 
\ee 
satisfying (\ref{4.5}) with 
$a=1/\overline{m}_{0}$, $b=1/\overline{m}_{-1}$ and $q=0=r$. 
Furthermore, 
the sequences $\{\lambda_i\}$ and $\{\kappa_i\}$ interlace: 
\be 
\lambda_1<\kappa_1<\lambda_2<
\cdots<\kappa_{l-1}<\lambda_{l}<\kappa_{l}<\lambda_{l+1}. 
                                             \label{5.8} 
\ee 
Conversely, for given $a,b>0$ and 
a discrete measure (\ref{5.7}), 
there exists a unique measure $m$ of the form (\ref{5.6}) 
satisfying (\ref{4.5}) with $q=0=r$. 
\\
(ii) Let $q>0$ be given. For every discrete measure 
\be 
m=\gamma_1\epsilon_{\lambda_1}
+\cdots+\gamma_{l}\epsilon_{\lambda_{l}} 
\quad \mbox{with} \quad  
\gamma_i>0 \ \mbox{and} \ 0<\lambda_1<\cdots<\lambda_{l}, 
                                             \label{5.9} 
\ee 
there exists a unique measure $M$ of the form (\ref{5.7}) 
satisfying (\ref{4.5}) with 
$r=0=a$ and $b=1/(q+\overline{m}_{-1})$. 
Furthermore, it holds that 
\be 
\lambda_1<\kappa_1<\lambda_2<
\cdots<\kappa_{l-1}<\lambda_{l}<\kappa_l
\le\lambda_l+\overline{m}_0 q^{-1}. 
                                             \label{5.10} 
\ee 
Conversely, for given $b>0$ and 
a discrete measure (\ref{5.7}), 
there exists a unique measure $m$ of the form (\ref{5.9}) 
satisfying (\ref{4.5}) with $q=1/(b+\overline{M}_0)$ and $r=0=a$. 
\end{pr}
{\it Proof.}~ 
(i) Observe from (\ref{4.5}) with $q=0$ that 
$g(\lambda):=\int(\lambda+u)^{-1}M(du)$ 
is a rational function of the form 
\[ 
g(\lambda)
=\frac{1}{\lambda}
\left(\frac{P(\lambda)}{Q(\lambda)}
-(a\lambda+b)\right)
=\frac{P(\lambda)-(a\lambda+b)Q(\lambda)}{\lambda}
\cdot\frac{1}{Q(\lambda)}, 
\] 
where $a$ and $b$ are given by (\ref{4.6}) and 
\[ 
P(\lambda)=\prod_{i=1}^{l+1}(\lambda+\lambda_i), 
\qquad 
Q(\lambda)=\sum_{i=1}^l\gamma_i\prod_{j\ne i}(\lambda+\lambda_j). 
\] 
We see also that 
$P_0(\lambda):=(P(\lambda)-(a\lambda+b)Q(\lambda))/\lambda$ 
is in fact a polynomial with degree less than or equal to $l-1$, and that 
$Q$ has a zero in the interval $(-\lambda_{i+1},-\lambda_i)$ 
for each $i\in\{1,\ldots,l\}$ because 
$Q(-\lambda_i)/Q(-\lambda_{i+1})<0$ as observed from (\ref{5.6}). 
Therefore, 
$g(\lambda)
=aP_0(\lambda)/((\lambda+\kappa_1)\cdots(\lambda+\kappa_{l}))$  
for some $\kappa_1,\ldots,\kappa_{l}$ satisfying (\ref{5.8}). 
It remains to find $c_1,\ldots,c_{l}>0$ such that 
\[ 
aP_0(\lambda)=\sum_{i=1}^{l}c_i\prod_{j\ne i}(\lambda+\kappa_j), 
\qquad \lambda>0.  
\] 
Since $\kappa_1,\ldots,\kappa_{l}$ are mutually different, 
a necessary and sufficient condition for 
the above identity to hold is that 
\be  
aP_0(-\kappa_i)=c_i\prod_{j\ne i}(-\kappa_i+\kappa_j), 
\qquad i\in\{1,\ldots,l\}, 
                                             \label{5.11} 
\ee 
which uniquely determine $c_1,\ldots,c_{l}$. 
Noting that $P_0(-\kappa_i)=P(-\kappa_i)/(-\kappa_i)$ 
because of $Q(-\kappa_i)=0$, 
we can verify the positivity of $c_i$'s 
by making use of (\ref{5.11}) and (\ref{5.8}).

It is almost a routine matter to show the converse assertion, 
whose proof we shall sketch. Let $M$ be given by (\ref{5.7}). 
Defining this time 
\[ 
P(\lambda)=\prod_{i=1}^{l}(\lambda+\kappa_i) 
\qquad \mbox{and} \qquad 
Q(\lambda)=\sum_{i=1}^{l}c_i\prod_{j\ne i}(\lambda+\kappa_j), 
\] 
we only have to show that the rational function 
\[ 
\left((a\lambda+b)+\lambda\frac{Q(\lambda)}{P(\lambda)}\right)^{-1}
           =:\frac{P(\lambda)}{Q_0(\lambda)}
\] 
can be rewritten into 
$\sum_{i=1}^{l+1}\gamma_i/(\lambda+\lambda_i)$ 
for some $\gamma_i>0$ and $\lambda_i>0$ satisfying (\ref{5.8}). 
Such $\lambda_2\ldots,\lambda_{l}$ are found as zeros of $Q_0$ 
since $Q_0(-\kappa_i)/Q_0(-\kappa_{i+1})<0$ 
\ $(i=1,\ldots,l-1)$. Therefore 
\be 
Q_0(\lambda)
=(\lambda+\lambda_2)\cdots(\lambda+\lambda_{l})Q_1(\lambda) 
                                             \label{5.12} 
\ee 
for some quadratic polynomial 
$Q_1(\lambda)=a\lambda^2+q_1\lambda+q_2$ with $q_1\in\R$ and $q_2>0$. 
Moreover, with the help of (\ref{5.12}) and (\ref{5.8}), 
we can show that 
$Q_1(-\kappa_i)<0$ for each $i\in\{1,\ldots,l\}$. 
These observations imply that 
$Q_1(\lambda)=a(\lambda+\lambda_1)(\lambda+\lambda_l)$ 
for some $\lambda_1\in (0,\kappa_1)$ and 
$\lambda_{l+1}\in (\kappa_{l},\infty)$. 
The rest of the proof (i.e. finding $\gamma_i$'s) 
is the same as before 
and the details are left to the reader. 

The assertion (ii) can be shown in an analogous way 
to the assertion (i). So we omit the proof of (ii) 
except the last inequality in (\ref{5.10}), 
which follows immediately from Lemma 4.2 (iii). 
\qed 

\subsection{Absolutely continuous Thorin measures}

\noindent 
As for continuous $m$, 
combining (\ref{4.8}) with the following inversion formula 
of Stieltjes-Perron (cf. \cite{Widder}, p.340, Corollary 7a) 
may provide explicit information on $M$: 
\[ 
M((s,t))+\frac{M(\{s\})+M(\{t\})}{2}
= 
-\frac{1}{\pi}\lim_{y\downarrow 0}\int_s^t{\rm Im}~\! G_M(x+{\rm i}y)dx, 
\quad -\infty<s<t<\infty, 
\] 
where ${\rm i}=\sqrt{-1}$. 
In particular, $M$ has a density function given by 
\[ 
\frac{dM}{dx}=-\frac{1}{\pi}\lim_{y\downarrow 0}{\rm Im}~\! G_M(x+{\rm i}y), 
\qquad x\in [s,t], 
\]  
provided that the right side converges boundedly and pointwise 
on $[s,t]$. 
We try to evaluate the right side in our setting. 
Let $m\in\cM$ be non-zero and set $q=0$ for simplicity. 
Rewrite (\ref{4.8}) in the form $G_M(z)=-1/(zG_m(z))+a-b/z$, 
where $a,b\ge 0$ and $M\in\cM$ correspond to 
the pair $(0,m)$ in the sense of Theorem 4.3. 
In the case where $m$ is absolutely continuous, 
it would be possible to establish a density formula for $M$ 
under suitable conditions, 
which do not seem, however, to be prescribed in a neat fashion. 
In typical cases, there exists a signed measure $m'$ 
on $[0,\infty)$ such that $G_m(z)=\int_{[0,\infty)}\log(z-u)m'(du)$ 
and hence 
\[ 
\lim_{y\downarrow 0}G_m(x+{\rm i}y)
=
\int_{[0,\infty)}\log|x-u|m'(du) 
+{\rm i}\pi m'((x,\infty)) 
\] 
for each $x>0$ with $m'(\{x\})=0$.
Therefore, the density formula would take the form 
\begin{eqnarray}
\frac{dM}{dx} 
& = & 
\frac{1}{\pi x}\lim_{y\downarrow 0}
{\rm Im}~\frac{1}{G_m(x+{\rm i}y)}              \nonumber \\ 
& = & 
-\frac{1}{x}\cdot
\frac{m'((x,\infty))}
{\ds{\left(\int_{[0,\infty)}\log|x-u|m'(du)\right)^2
+\left(\pi m'((x,\infty))\right)^2}} 
                                             \label{5.13} 
\end{eqnarray} 
for $x>0$. Similarly, (\ref{4.8}) in principle makes it possible 
to derive the information of $m$ 
corresponding to given $a,b\ge 0$ and $M\in\cM$. 
We shall give some concrete examples of such a kind, 
in which we continue to take $q=0$. 

\bigskip 

\noindent 
{\it Example 5.3}~(i) 
For $0<\alpha<1$ and $\kappa\ge 0$, 
let $m$ be as in (\ref{4.16}). 
Then $G_m(z)=-(\kappa-z)^{-\alpha}$. 
Here the power function $z^{\alpha}$ 
of $z\in\C\setminus (-\infty,0]$ is defined 
to be $|z|^{\alpha}(\cos \arg z+{\rm i}\sin \arg z)$ 
with $\arg z$ chosen so that $|\arg z|<\pi$. 
Since we know from Lemma 4.2 (iii) that 
$\cS(M)\subset [\kappa, \infty)$, 
fix an $x>\kappa$ arbitrarily. 
It is easy to see that 
\[ 
\frac{1}{\pi x}\lim_{y\downarrow 0}
{\rm Im}~\frac{1}{G_m(x+{\rm i}y)}
= 
-\frac{1}{\pi x}(x-\kappa)^{\alpha}
\sin(\alpha(-\pi))
= 
\frac{1}{x}\cdot\frac{(x-\kappa)^{\alpha}}
{\Gamma(\alpha)\Gamma(1-\alpha)}. 
\] 
It is not difficult to verify that 
the above convergence is uniform in $x$ 
on every compact interval contained in $(\kappa, \infty)$. 
Thus (\ref{4.17}) has been recovered. \\ 
(ii) 
Let $0\le \lambda_1<\lambda_2$ be arbitrary. 
Define $m(du)=\one_{(\lambda_1, \lambda_2)}(u)du$, 
for which $a=1/(\lambda_1-\lambda_2)$, 
$b=1/(\log\lambda_2-\log\lambda_1)$. Moreover, 
$G_m(z)=\log(z-\lambda_1)-\log(z-\lambda_2)$ and 
hence $m'=\epsilon_{\lambda_1}-\epsilon_{\lambda_2}$. 
By (\ref{5.13}) 
\be 
\frac{dM}{dx} 
= 
\frac{1}{x}\cdot
\frac{\one_{(\lambda_1, \lambda_2)}(x)}
{\left(\ds{\log\frac{x-\lambda_1}{\lambda_2-x}}\right)^2+\pi^2}. 
                                             \label{5.14} 
\ee 
While (\ref{4.7}) tells that 
\[ 
M([\lambda_1, \lambda_2])
=
\frac{\lambda_1+\lambda_2}{2(\lambda_2-\lambda_1)}
-\frac{1}{\log\lambda_2-\log\lambda_1}, 
\] 
it is not clear how to verify this 
directly from (\ref{5.14}) unless $\lambda_1=0$. 
If $\lambda_1>0$, we have also 
the sector constant estimate (\ref{4.11}) 
which reads 
${\rm Sect}(\cE^{\delta})-1 \le 
\sqrt{(\lambda_2-\lambda_1)/\lambda_1}$ 
for each $\delta>0$. 
For the special choice $\lambda_1=0$ and $\lambda_2=1$, 
the density function of the GGC with pair $(0,m)$ 
can be found in \cite{JRY08} (Eq.(259)). \\  
(iii) 
Let $a,b>0$ be given arbitrarily. 
Consider $M(du)=\one_{(0,1)}(u)du$ for simplicity.  
The density of the absolutely continuous part $m_c$ of 
the Thorin measure $m$ corresponding to $(a,b,M)$  
in the sense of Theorem 4.4 (ii) 
can be calculated from (\ref{4.8}) with $q=0=r$, 
namely $G_m(z)=1/(az-b-zG_M(z))$, 
where $G_M(z)=\log z-\log(z-1)$. 
Indeed, defining $H(z)=a-(b/z)-G_M(z)$ 
for $z\ne 0$, one can show that, for each $x>0$ with $x\ne 1$, 
$H(x+{\rm i}y) \to a-(b/x)-\log|x/(x-1)|+{\rm i}\pi \one_{(0,1)}(x)$ 
as $y\downarrow 0$. Therefore, 
\[ 
\frac{dm_c}{dx}
=
-\frac{1}{\pi x}\lim_{y\downarrow 0}
{\rm Im}~\! \frac{1}{H(x+{\rm i}y)} 
= 
\frac{1}{x}\cdot
\frac{\one_{(0,1)}(x)}
{\left(\ds{a-\frac{b}{x}-\log\frac{x}{1-x}}\right)^2+\pi^2}.  
\] 
The point which requires extra care is the unique pole, 
say $x_0$, of $G_m$, which is located on 
the interval $(1,\infty)$. 
It is characterized as a unique solution to 
\be 
ax-b+x\log\left(1-\frac{1}{x}\right)=0, 
                                             \label{5.15} 
\ee from which $b/a<x_0<1+(b+1)/a$ 
can be deduced with the help of elementary inequalities 
$u/(1+u)<\log(1+u)<u$ for $u>-1$. 
The point mass of $m$ at $x_0$ is given as 
the residue of $G_m$ at $z=x_0$: 
\begin{eqnarray*} 
m(\{x_0\})
& = & 
\lim_{z\to x_0}G_m(z)(z-x_0)
\ = \ 
\left(
\left.\frac{d}{dz}\left(az-b-zG_M(z)\right)\right|_{z=x_0}
\right)^{-1} \\ 
& = &  
\left(a+\log\left(1-\frac{1}{x_0}\right)+\frac{1}{x_0}\right)^{-1}
\ = \ 
\left(\frac{b}{x_0}+\frac{1}{x_0-1}\right)^{-1}, 
\end{eqnarray*} 
where (\ref{5.15}) has been used to get the last expression. 

\subsection{Connections with noncommutative probability theory} 
In the previous subsection, 
the calculus of Stieltjes transforms played 
quite an important role. 
So, it might be no surprise that 
observations we had made so far 
have some connections with 
noncommutative probability theory, 
another context in which the reciprocal of 
the Stieltjes transform 
serves as one of essential tools. 
We will be particularly concerned with 
the Boolean convolution and the free Poisson distributions 
(known also as the Marchenko-Pasture laws). 
To give the definition of the former, 
we follow \cite{SW} and 
introduce the `Boolean cumulant' 
$K_m(z):=z-1/G_{m}(z)$ for a {\it probability} measure $m$ on $\R$. 
For our purpose the domain of such an operation 
shall be restricted to $\cM_1$, 
the totality of probability measures on $(0,\infty)$. 
For any $m_1,m_2\in\cM_1$, 
the Boolean convolution $m_1\uplus m_2$ of $m_1$ and $m_2$ 
is then defined to be an element of $\cM_1$ such that 
\be 
K_{m_1\uplus m_2}(z)=K_{m_1}(z)+K_{m_2}(z), 
\qquad z\in\C\setminus\R_+.  
                                             \label{5.16} 
\ee 
Let $t>0$ be arbitrary. Following \cite{BW}, 
one can define also the $t$th Boolean convolution power 
$m^{\uplus t}$ of $m\in\cM_1$ by 
\be 
K_{m^{\uplus t}}(z)=t \cdot K_{m}(z), 
\qquad z\in\C\setminus\R_+.  
                                             \label{5.17} 
\ee 
It is a good exercise to verify from Lemma 4.1 with $a=1$ 
that the requirements (\ref{5.16}) and (\ref{5.17}) determine uniquely 
such measures $m_1\uplus m_2$ and $m^{\uplus t}$, 
respectively. (In fact, for the verification, 
one needs to observe in Lemma 4.1 additionally that 
$\overline{M}_{-1}=\infty$ whenever $b=0$. 
But this can be seen from (\ref{4.25}).) 
The following proposition, which may have a number of variants, 
can be regarded essentially 
as a reformulation of this fact 
in the language of GGC's and the corresponding CBCI processes. 
\begin{pr}
Suppose that probability measures $m_1$ and $m_2$ on $(0,\infty)$ 
are Thorin measures. For each $i\in\{1,2\}$, 
let the quadruplet $(1,b_i,n_i,1)$ is the one 
determined from the pair $(0,m_i)$ by Theorem 4.3. 
Then the following assertions hold. \\ 
(i) $m_1\uplus m_2$ is a Thorin measure 
and the GGC with pair $(0,m_1\uplus m_2)$ 
is a unique 
stationary distribution of the CBCI-process 
with quadruplet $(1,b_1+b_2,n_1+n_2,1)$. \\ 
(ii) For each $t>0$, $m_1^{\uplus t}$ is a Thorin measure 
and the GGC with pair $(0,m_1^{\uplus t})$ 
is a unique 
stationary distribution of the CBCI-process 
with quadruplet $(1,tb_1,t n_1,1)$. 
\end{pr}
{\it Proof.}~ 
(i) According to (\ref{4.10}), $n_i$ are of the form 
$n_i(dy)=dy\int u^2e^{-uy}M_i(dy)$ for some $M_i\in\cM$. 
Observe from (\ref{4.5}) with $a=1$ that 
\be 
\left(\int\frac{m_i(du)}{\lambda+u}\right)^{-1}-\lambda
=b_i+\lambda\int\frac{M_i(du)}{\lambda+u}\ge 0, 
\qquad \lambda>0. 
                                             \label{5.18} 
\ee 
Combining the above inequality with (\ref{5.16}), we get 
\[ 
\left(\int\frac{(m_1\uplus m_2)(du)}{\lambda+u}\right)^{-1}-\lambda
\ge \left(\int\frac{m_1(du)}{\lambda+u}\right)^{-1}-\lambda, 
\qquad \lambda>0 
\] 
or 
$\int(\lambda+u)^{-1}(m_1\uplus m_2)(du)\le 
\int(\lambda+u)^{-1}m_1(du)$ for $\lambda>0$. 
By integrating with respect to 
the Lebesgue measure $d\lambda$ over $[0,1]$ 
and then applying Fubini's theorem 
$\int\log(1+u^{-1})(m_1\uplus m_2)(du)
\le \int\log(1+u^{-1})m_1(du)<\infty$. 
Hence $m_1\uplus m_2$ is a Thorin measure. 
The last half of the assertion follows from 
the equalities for $i=1$ and $i=2$ in (\ref{5.18}). 
Indeed, summing up them leads to 
\[ 
\left(\int\frac{(m_1\uplus m_2)(du)}{\lambda+u}\right)^{-1}-\lambda
=(b_1+b_2)+\lambda\int\frac{(M_1+M_2)(du)}{\lambda+u} 
\] 
or 
\[ 
\int\frac{(m_1\uplus m_2)(du)}{\lambda+u}
=
\frac{1}{\lambda+(b_1+b_2)+\ds{\lambda\int\frac{(M_1+M_2)(du)}{\lambda+u}}} 
\] 
for all $\lambda>0$. 
This is sufficient for the proof of (i). 

The proof of (ii) proceeds along the same lines 
as that of (i) on noting that by (\ref{5.17}) and (\ref{5.18}) 
\[ 
\int\frac{m_1^{\uplus t}(du)}{\lambda+u} 
= 
\frac{1}{\lambda+tb_1+\ds{\lambda t\int\frac{M_1(du)}{\lambda+u}}} 
\le 
\frac{\max\{t^{-1},1\}}
{\lambda+b_1+\ds{\lambda \int\frac{M_1(du)}{\lambda+u}}}  
\] 
for any $\lambda>0$. The details are omitted. 
\qed 

\bigskip 

The final topic is related to the free probability theory. 
In that theory, the counterpart (in an appropriate sense) 
of the Poisson distribution exists and 
is called the free Poisson distribution. 
It is, by definition, of the form 
\[ 
P_{\alpha,\beta}(du)
:= 
\left\{
\begin{array}{ll}
(1-\beta)\epsilon_{0}(du)+\beta p_{\alpha,\beta}(u)du 
&  (0\le \beta < 1) \\ 
p_{\alpha,\beta}(u)du 
&  (\beta\ge 1)  
\end{array}
\right.
\] 
for some $\alpha>0$ and $\beta\ge 0$, where 
\[ 
p_{\alpha,\beta}(u)
= 
\frac{1}{2\pi\alpha u}
\sqrt{4\alpha^2\beta-\left(u-\alpha(1+\beta)\right)^2}
\one_{[\alpha(1-\sqrt{\beta})^2,\alpha(1+\sqrt{\beta})^2]}(u). 
\] 
Notice  that $P_{\alpha,\beta}$ is a Thorin measure 
if and only if $\beta\ge 1$. 
The formula for the Stieltjes transform of this distribution is 
\be 
G_{P_{\alpha,\beta}}(z)
= 
\frac{\ds{z+\alpha(1-\beta)-
\sqrt{(z+\alpha(1-\beta))^2-4\alpha z}}}{2\alpha z}. 
                                           \label{5.19} 
\ee 
(See e.g. p.205 in \cite{NS}.) 
We now remark that a class of the free Poisson distributions plays 
a special role in describing the fixed points of 
the correspondence between $m$ and $M$ defined through Lemma 4.1 
although we don't have any interpretation 
in the context of CBI-processes. 
Here is an explicit statement. 
\begin{pr}
Let $q\ge 0$ and suppose that $m\in\cM$ is non-zero. 
Assume that $a,b\ge 0$ and $M\in\cM$ satisfy 
(\ref{4.5}) with $r=0$. Then $m$ coincides with $M$ if and only if 
$q=0$ or $\overline{m}_0<\infty$ and $m$ is given by 
\be 
m(du)=
\left\{
\begin{array}{ll} 
\ds{\frac{1-bq}{a+q}}P_{\alpha,\beta}(du) 
& (\overline{m}_0<\infty) \\ 
\ds{\frac{1}{\pi u}}\sqrt{u-(b/2)^2}\one_{[(b/2)^2,\infty)}(u)du 
& (q=0, \ \overline{m}_0=\infty), 
\end{array}
\right. 
                                           \label{5.20}
\ee 
where $\alpha=(1-bq)/(a+q)^2$ and $\beta=(1+ab)/(1-bq)$. 
\end{pr}
{\it Proof.}~ 
In view of (\ref{4.8}), it is obvious that 
$m=M$ if and only if $G_m(z)-q=1/(az-b-zG_m(z))$. 
By solving it we can deduce 
\[ 
G_m(z)
=
\frac{(a+q)z-b-\sqrt{((a+q)z-b)^2-4(1-bq)z}}{2z}. 
\] 
In the case where $\overline{m}_0<\infty$, 
we have $a+q>0$ by (\ref{4.6}), 
and the proof concludes by comparing this with (\ref{5.19}). 
If $\overline{m}_0=\infty$ (and hence $\overline{M}_0=\infty$), 
then (\ref{4.7}) and (\ref{4.6}) imply $q=0$ and $a=0$, respectively. 
Consequently, we have $G_m(z)=(-b-\sqrt{b^2-4z})/(2z)$, 
from which the second expression in (\ref{5.20}) 
can be derived by inversion. \qed 

\bigskip 

\noindent 
Denoting by $\rho_{q,a,b}(u)$ and $\rho_b(u)$ the densities 
of the measures on the right side of (\ref{5.20}) 
in the first and the second cases, respectively, 
we note that 
\[ 
\lim_{a\downarrow 0}\rho_{0,a,b}(u)
=\rho_{b}(u)
=\lim_{q\downarrow 0}\rho_{q,0,b}(u) 
\] 
for each $u>0$ and $b\ge 0$. 
Remark also that $\rho_0(u)du$ coincides with $m(du)$ 
in (\ref{4.16}) with $\alpha=1/2$ and $\kappa=0$, 
for which clearly $m=M$ holds. 

\bigskip 

\noindent 
{\bf Acknowledgment.}~ 
The author is grateful to Professors M. Yor and H. Masuda 
for bringing to his attention \cite{SSV} and \cite{KRM}, respectively.

\end{document}